\documentclass[11pt]{article}
\usepackage[authoryear]{natbib}
\newtheorem{theorem}{Theorem}
\newtheorem{lemma}{Lemma}
\newtheorem{proposition}{Proposition}
\newtheorem{corollary}{Corollary}
\newtheorem{definition}{Definition}

\newtheorem{assumption}{Assumption}

\newcommand{\lab}[1]{\label{#1}}
\newcommand{\labs}[1]{\label{#1}}
\newcommand{\labe}[1]{\label{#1}}

\newcommand{\dating}[3]{\date{{\Large {\bf #1}}\vspace{5mm}\\{\Large#3}}\label{#2}}

\usepackage{amssymb}
\def\Bbb{\mathbb}

\newcommand{\Section}[1]{\section{#1}\setcounter{equation}{0}}

\newcommand{\cqfd}{\qquad\framebox[2.7mm]{\rule{0mm}{.7mm}}}
\newcommand{\bm}[1]{\mbox{\boldmath $#1$}}
\newcommand{\scal}[2]{\langle #1,#2\rangle}

\newcommand{\st}{\strut}
\newcommand{\1}{1\hskip-2.6pt{\rm l}}

\def\tr{{\rm tr}}
\def\gbeta{{\bm{\beta}}}
\def\gq{{\bm q}}
\def\galpha{{\bm{\alpha}}}
\def\gp{{\bf p}}
\def\gR{{\bf R}}
\def\gL{{\bf L}}
\def\RR{{\mathcal R}}

\def\gtheta{{\bm{\theta}}}
\def\w{{\rm w}}

\def\F{{\mathcal F}}
\def\E{{\mathbb E}}

\def\H{{\mathcal H}}
\def\B{{\mathcal B}}
\def\P{{\mathbb P}}
\def\PP{{\mathcal P}}
\def\N{{\mathbb{N}}}
\def\IL{{\Bbb{L}}}
\def\R{{\mathbb R}}
\def\S{{\mathcal S}}

\def\SS{{\mathbb S}}
\def\BB{\mathbb B}

\def\HH{{\mathbb H}}
\def\H{{\mathcal H}}
\def\TT{{\mathbb T}}
\def\gT{{\bf T}}
\def\T{{\mathcal T}}

\def\FF{{\mathbb F}}
\def\C{\mathcal C}
\def\D{\mathcal D}

\def\L{{\mathcal L}}
\def\risq{{\mathcal R}}
\def\gR{{\bm R}}

\newcommand{\norm}[1]{\left\|{#1}\right\|}
\def\1{1\hskip-2.6pt{\rm l}}
\def\eps{{\varepsilon}}
\def\<{{\langle}}
\def\>{{\rangle}}

\newcommand{\eref}[1]{(\ref{#1})}
\newcommand{\pa}[1]{\left({#1}\right)}
\newcommand{\cro}[1]{\left[{#1}\right]}
\newcommand{\ab}[1]{\left|{#1}\right|}
\newcommand{\ac}[1]{\left\{{#1}\right\}}

\begin{document}

\title{\mbox{}\vspace{-10mm}\\{\Huge{\bf  Estimating\vspace{2mm} composite functions
by model selection\vspace{5mm}}}}
\author{
{\Large Yannick Baraud\vspace{3mm}}\\{\Large Universit\'e de Nice Sophia-Antipolis
\vspace{1mm}}\\{\Large Laboratoire J-A Dieudonn\'e}\vspace{5mm}\\
{\Large Lucien Birg\'{e}\vspace{3mm}}\\{\Large Universit\'e Paris VI\vspace{1mm}}\\
{\Large Laboratoire de Probabilit\'es et Mod\`eles Al\'eatoires\vspace{1mm}}\\
{\Large U.M.R.\ C.N.R.S.\ 7599}}\vspace{5mm}
\dating{\mbox{}}{June 2006}{February  9, 2011}      
\maketitle 

\footnotetext{{\it AMS 1991 subject classifications.} Primary 62G05

{\it \hspace{2.3mm}Key words and phrases.} Curve estimation, model selection, composite 
functions.}

\begin{abstract}
We consider the problem of estimating a function $s$ on $[-1,1]^{k}$ for large values of $k$ by 
looking for some best approximation by composite functions of the form $g\circ u$. Our solution is
based on model selection and leads to a very general approach to solve this problem with respect too
many different types of functions $g,u$ and statistical frameworks. In particular, we handle the
problems of approximating $s$ by additive functions, single and multiple index models, neural
networks, mixtures of Gaussian densities (when $s$ is a density) among other examples. We also
investigate the situation where $s=g\circ u$ for functions $g$ and $u$ belonging to possibly
anisotropic smoothness classes. In this case, our approach leads to a completely adaptive estimator
with respect to the regularity of $s$. 
\end{abstract}

\Section{Introduction\labs{I}}
In various statistical problems, we have at hand a random mapping $\bm{X}$ from a
measurable space $(\Omega,{\cal A})$ to $(\Bbb{X},{\cal X})$ with an unknown
distribution $P_s$ on $\Bbb{ X}$ depending on some parameter $s\in{\cal S}$ which is a function 
from $[-1,1]^k$ to $\Bbb{R}$. For instance, $s$ may be the density of an i.i.d.\ sample or the
intensity of a Poisson process on  $[-1,1]^k$ or a regression function. The statistical problem amounts
to estimating $s$ by some estimator $\widehat{s}=\widehat{s}(\bm{X})$ the performance of which
is measured by its quadratic risk,
\[
R(s,\widehat{s})=\Bbb{E}_s\!\left[d^2(s,\widehat{s})\right],
\]
where $d$ denotes a given distance on ${\cal S}$. To be more specific,  we shall assume in this
introduction that $\bm{X}=(X_{1},\ldots,X_{n})$ is a sample of density $s^{2}$ (with $s\ge0$) with
respect to some measure $\mu$ and $d$ is the Hellinger distance. We recall that, given two
probabilities $P,Q$ dominated by $\mu$ with respective densities $f=dP/d\mu$ and
$g=dQ/d\mu$, the Hellinger distance $h$ between $P$ and $Q$ or, equivalently, between
$f$ and $g$ (since it is independant of the choice of $\mu$ is given by 
\begin{equation}
h^2(P,Q)=h^2(f,g)=\frac{1}{2}\int\left(\sqrt{f}-\sqrt{g}\right)^2d\mu.
\labe{Eq-Hell}
\end{equation}
It follows that  $\sqrt{2}d(s,t)$ is merely the $\IL_{2}$-distance between $s$ and $t$.  A
general method for constructing estimators $\widehat{s}$ is to choose a model $S$ for $s$, that is,
do as if $s$ belonged to $S$ and to build $\widehat{s}$ as an element of $S$. Sometimes the
statistician really assumes that $s$ belongs to $S$ and that $S$ is the true parameter set,
sometimes he does not and rather considers
$S$  as an approximate model. This latter approach is somewhat more reasonable since it is in
general impossible to be sure that $s$ does belong to $S$. Given $S$ and a suitable estimator
$\widehat{s}$, as those built in Birg\'e~\citeyearpar{MR2219712} for example, one can achieve a
risk bound of the form
\begin{equation}
R(s,\widehat{s})\le C\left[\inf_{t\in S}d^2(s,t)+\tau\D(S)\right],
\labe{Eq-RB1}
\end{equation}
where $C$ is a universal constant (independent of $s$), $\D(S)$ the dimension of the model
$S$ (with a proper definition of the dimension) and $\tau$, which is equal to $1/n$ in the 
specific context of density estimation,  characterizes the amount of information provided by the
observation $\bm{X}$. 

It is well-known that many classical estimation procedures suffer from the so-called
``curse of dimensionality", which means that the risk bound (\ref{Eq-RB1}) deteriorates
when $k$ increases and actually becomes very loose for even moderate values of $k$. This
phenomenon is easy to explain and actually connected with the most classical way of choosing 
models for $s$. Typically, and although there is no way to check that such an assumption is 
true, one assumes that $s$ belongs to some smoothness class (H\"older, Sobolev or Besov) of
index $\alpha$ and such an assumption can be translated in terms of approximation properties
with respect to the target function $s$ of a suitable collection of linear spaces (generated by
piecewise polynomials, splines, or wavelets for example).   More precisely, there exists a
collection $\SS$ of models  with the following property: for all $D\ge 1$, there exists a model
$S\in\SS$ with dimension $D$ which approximates $s$ with an error bounded by
$cD^{-\alpha/k}$ for some $c$ independent of $D$ (but depending on $s$, $\alpha$ and $k$). 
With such a collection at hand, we deduce from~\eref{Eq-RB1} that whatever $D\ge 1$ one can 
choose a model $S=S(D)\in\SS$ for which the estimator $\widehat s\in S$ achieves a risk 
bounded from above by $C\left[c^2D^{-2\alpha/k}+\tau D\right]$. Besides, by using the 
elementary Lemma~\ref{L-opt} below to be proved in Section~\ref{H5}, one can optimize the 
choice of $D$, and hence of the model $S$ in $\SS$, to build an estimator whose risk satisfies
\begin{equation}
R(s,\widehat{s})\le
C\max\left\{3c^{2k/(2\alpha+k)}\tau^{2\alpha/(2\alpha+k)};2\tau\right\}.
\labe{eq00}
\end{equation}
%
\begin{lemma}\lab{L-opt}
For all positive numbers $a$, $b$ and $\theta$ and $\Bbb{N}^\star$ the set of  positive integers,
\[
\inf_{D\in\Bbb{N}^\star}\{aD^{-\theta}+bD\}\le
b+\min\left\{2a^{1/(\theta+1)}b^{\theta/(\theta+1)};a\right\}\le
\max\left\{3a^{1/(\theta+1)}b^{\theta/(\theta+1)};2b\right\}.
\]
\end{lemma}
Since the risk bound~\eref{eq00} is achieved for $D$ of order $\tau^{-k/(2\alpha+k)}$, as 
$\tau$ tends to 0, the deterioration of the rate $\tau^{2\alpha/(2\alpha+k)}$ when $k$
increases comes from the fact that we use  models of larger dimension to approximate $s$ when
$k$ is large. Nevertheless,  this phenomenon is only  due to the previous approach based on
smoothness assumptions  for $s$. An alternative approach, assuming that $s$ can be closely
approximated by suitable parametric models the dimensions of which do not depend on $k$
would not suffer from the same weaknesses. More generally, a structural assumption on $s$
associated to a collection of models $\SS'$, the approximation properties of which improve on 
those of $\SS$, can only lead to a better risk bound and it is not clear at all that assuming that 
$s$ belongs to a smoothness class is more realistic than directly assuming approximation 
bounds with respect to the models of $\SS'$. Such structural assumptions that would amount to
replacing the large models involved in the approximation of smooth functions by simpler ones
have been used for many years, especially in the context of regression. Examples of such
structural assumptions are provided by additive models, the single index model, models based
on radial approximation as in Donoho and Johnstone~\citeyearpar{MR981438}, the projection 
pursuit algorithm  indroduced by Friedman and Tuckey~\citeyearpar{Freidman74}, (an
overview of the procedure is available in Huber~\citeyearpar{MR790553}), and artificial neural
networks as in Barron~\citeyearpar{MR1237720,Barron1994}, among other examples. 

In any case, an unattractive feature of the previous approach based on an a priori choice of a
model $S\in\SS$ is that it requires to know suitable  upper bounds on the distances between $s$
and the models $S$ in $\SS$. Such a requirement is much too strong and an essential
improvement can be brought by the modern theory of model selection. More precisely, given
some prior probability $\pi$ on $\SS$, model selection allows to build  an estimator
$\widehat{s}$ with a risk bound
\begin{equation}
R(s,\widehat{s})\le C\inf_{S\in\SS}\ac{\inf_{t\in S}d^2(s,t)+
\tau\cro{\D(S)+\log\pa{1/\pi(S)\st}}}, 
\labe{Eq-RB2}
\end{equation}
for some universal constant $C$. If we neglect the influence of $\log\pa{1/\pi(S)}$, which is
connected to the complexity of the family $\SS$ of models we use, the comparison between
(\ref{Eq-RB1}) and (\ref{Eq-RB2}) indicates that the method approximately selects among 
$\SS$ a model leading to the smallest risk bound.

With such a tool at hand, that allows us to play with many models simultaneously and let the
estimator choose a suitable one, we may freely introduce various models corresponding to
various sorts of structural assumptions on $s$ that avoid the ``curse of dimensionality''. We can,
moreover, mix them with models that are based on pure smoothness assumptions that do suffer
from this dimensional effect or even with simple parametric models. 

The main purpose of this paper is to provide a method for building various sorts of
models that may be used, in conjonction with other ones, to approximate functions on
$[-1,1]^k$ for large values of $k$. The idea, which is not new, is to approximate the
unknown $s$ by a composite function $g\circ u$ where $g$ and $u$ have different
approximation properties. Recent works in this direction can be found in Horowitz and 
Mammen~\citeyearpar{MR2382659} or Juditsky, Lepski and Tsybakov~\citeyearpar{MR2509077}. 
Actually, our initial motivation for this research was a series of lectures given at CIRM in 2005 by
Oleg Lepski about a former version of this last  paper. There are, nevertheless, major differences
between their approach and ours. They deal with estimation in the white noise model, kernel
methods and the $\Bbb{L}_\infty$-loss. They also assume that the true unknown density $s$ to
be estimated can be written as $s=g\circ u$ where $g$ and $u$ have given smoothness
properties and use these properties to build a kernel estimator which is better than those based on
the overall smoothness of $s$. The use of the $\Bbb{L}_\infty$-loss indeed involves
additional difficulties and the minimax rates of convergence happen to be substantially slower
(not only by logarithmic terms) than the rates one gets for the $\Bbb{L}_2$-loss, as the authors
mention on page 1369, comparing their results with those of Horowitz and Mammen (2007).

Our approach is radically different from the one of  Juditsky, Lepski and Tsybakov and
considerably more general as we shall see, but this level of generality has a price. While they
provide a constructive estimator that can be computed in a reasonable amount of time, although
based on supposedly known smoothness properties of $g$ and $u$, we offer a general but
abstract method that applies to many situations but does not provide practical estimators, only
abstract ones. As a consequence, our results about the performance of these estimators are of a
theoretical nature, to serve as benchmarks about what can be expected from good estimators in
various situations. 

We actually consider ``curve estimation" with an unknown functional parameter $s$ and
measure the loss by $\Bbb{L}_2$-type distances. Our construction applies to various statistical
frameworks (not only the Gaussian White Noise but all these for which a suitable model
selection theorem is available). We also do not assume that $s=g\circ u$ but rather approximate
$s$ by functions of the form $g\circ u$ and do not fix in advance the smoothness properties of
$g$ and $u$ but rather let our estimator adapt to it. This approach leads to a completely adaptive
method with many different possibilities to approximate $s$. It allows, in particular, to play with
the smoothness properties of $g$ and $u$ or to mix purely parametric models with others based
on smooth functions. Since methods and theorems about model selection are already available,
our main task here will be to build suitable models for various forms of composite functions
$g\circ u$ and check that they do satisfy the assumptions required for applying previous model
selection results.

\Section{Our statistical framework\labs{F}}
We observe a random element $\bm{X}$ from the probability space $(\Omega,{\cal A},\P_{s})$ 
to $(\Bbb{X},{\cal X})$ with distribution $P_s$ on $\Bbb{X}$ depending on an unknown
parameter $s$. The set ${\cal S}$ of possible values of $s$ is a subset of some space
$\Bbb{L}_q(E,\mu)$ where $\mu$ is a given {\em probability} on the measurable space
$(E,{\cal E})$. We shall mainly consider the case $q=2$ even though one can  also take $q=1$ in
the context of density estimation.  We denote by $d$ the distance on $\Bbb{L}_q(E,\mu)$
corresponding to the $\Bbb{L}_q(E,\mu)$-norm $\|\cdot\|_q$ (omitting the dependency of
$d$ with respect to $q$) and by $\Bbb{E}_s$ the expectation with  respect to $\P_{s}$ so that the
quadratic risk of an estimator $\widehat{s}$ is $\Bbb{E}_s\!\left[d^2\left(s,\widehat{s}\right)
\st\right]$. The main objective of this paper, in order to estimate $s$ by model selection, is to
build special models $S$ that consist of functions of the form $f\circ t$ where $t=(t_1\ldots,t_l)$
is a mapping from $E$ to $I\subset\Bbb{R}^l$, $f$ is a continuous function on $I$ and
$I=\prod_{j=1}^lI_j$ is a product of compact intervals of $\Bbb{R}$. Without loss of generality,
we may assume that
$I=[-1,1]^l$. Indeed, if  $l=1$, $t$ takes its values in $I_1=[\beta-\alpha,\beta+\alpha]$,
$\alpha>0$ and $f$ is defined on $I_1$, we can replace the pair $(f,t)$ by $(\bar{f},\bar{t})$
where $\bar{t}(x)=\alpha^{-1} [t(x)-\beta]$ and $\bar{f}(y)=f(\alpha y+\beta)$ so that $\bar{t}$
takes its values in $[-1,1]$ and $f\circ t=\bar{f}\circ\bar{t}$. The argument easily extends to the
multidimensional case.

\subsection{Notations and conventions\labs{F5}}
To perform our construction based on composite functions $f\circ t$, we introduce the
following spaces of functions : $\T\subset\Bbb{L}_q(E,\mu)$ is the set of measurable
mappings from $E$ to $[-1,1]$, $\F_{l,\infty}$ is the set of  bounded functions on $[-1,1]^l$ 
endowed with the distance $d_\infty$ given by $d_\infty(f,g)=\sup_{x\in[-1,1]^l} |f(x)-g(x)|$
and $\F_{l,c}$ is the subset of $\F_{l,\infty}$ which consists of continuous functions on 
$[-1,1]^l$. We denote by $\Bbb{N}^\star$ (respectively, $\R_{+}^{\star}$) the set of  positive
integers (respectively positive numbers) and set
\[ 
\lfloor z\rfloor=\sup\{j\in\Bbb{Z}\,|\,j\le z\}\quad\mbox{and}\quad
\lceil z\rceil=\inf\{j\in\Bbb{N}^\star\,|\,j\ge z\},\;\mbox{ for all }z\in\Bbb{R}.
\] 
The numbers $x\wedge y$ and $x\vee y$ stand for $\min\{x,y\}$ and $\max\ac{x,y}$ 
respectively and $\log_+(x)$ stands for $(\log x)\vee0$. The cardinality of a set $A$ is denoted 
by $|A|$ and, by convention,  ``countable" means ``finite or countable". We call {\em 
subprobability} on some countable set $A$ any positive measure $\pi$ on $A$ with $\pi(A)\le1$ 
and, given $\pi$ and $a\in A$, we set $\pi(a)=\pi(\{a\})$ and $\Delta_{\pi}(a)=-\log(\pi(a))$ with the 
convention $\Delta_{\pi}(a)=+\infty$ if $\pi(a)=0$. The dimension  of the linear space $V$ is 
denoted by $\D(V)$. Given a compact subset $K$ of $\Bbb{R}^k$ with $\stackrel{\circ}{K}\ne
\emptyset$, we define the {\em Lebesgue probability } $\mu$ on $K$ by $\mu(A)=\lambda(A)/
\lambda(K)$ for $A\subset K$,  where $\lambda$ denotes  the Lebesgue measure on 
$\Bbb{R}^k$.

For $x\in\R^m$, $x_j$ denotes the $j^{\rm th}$ coordinate of $x$ ($1\le j\le m$) and, similarly, 
$x_{i,j}$ denotes the $j^{\rm th}$ coordinate of $x_i$ if the vectors $x_i$ are already indexed. We 
set $\ab{x}^{2}=\sum_{j=1}^mx_j^{2}$ for the squared Euclidean norm of $x\in\Bbb{R}^m$, 
without reference to the dimension $m$, and denote by ${\cal B}_m$ the corresponding unit ball 
in $\R^m$ and by $\bm{1}$ the vector with unit coordinates in $\Bbb{R}^m$. Similarly, 
$\ab{x}_{\infty}=\max\{|x_{1}|,\ldots,|x_{m}|\}$ for all $x\in\R^{m}$. For $x$ in some metric space
$(M,d)$ and $r>0$, ${\cal B}(x,r)$ denotes the closed ball of center $x$ and radius $r$ in $M$ 
and for $A\subset M$, $d(x,A)=\inf_{y\in A} d(x,y)$. Finally, $C$ stands for a universal constant 
while $C'$ is a constant that depends on some parameters of the problem. We may make this 
dependence explicit by writing $C'(a,b)$ for instance. Both $C$ and $C'$ are generic notations 
for constants that may change from line to line.

\subsection{A general model selection result\labs{F1}}
General model selection results apply to models which possess a {\em finite dimension} in a
suitable sense. Throughout the paper, we assume that in the statistical framework we
consider the following theorem holds.
%
\begin{theorem}\lab{T-models}
Let $\SS$ be a countable family of finite dimensional linear subspaces $S$ of
$\Bbb{L}_q(E,\mu)$ and let $\pi$ be some subprobability measure on $\SS$.  There exists an
estimator $\widehat{s}=\widehat{s}(\bm{X})$ with values in $\cup_{S\in\SS}S$ satisfying, for
all $s\in{\cal S}$,
\begin{equation}
\Bbb{E}_s\!\left[d^2\left(s,\widehat{s}\right)\st\right]\le C\inf_{S\in\SS}
\left\{d^2\left(s,S\right)+\tau\left[\st(\D(S)\vee1)+\Delta_{\pi}(S)\right]\right\},
\labe{Eq-main1}
\end{equation}
where the positive constant $C$ and parameter $\tau$ only depend on the specific statistical 
framework at hand.
\end{theorem}
%
Similar results often hold also for the loss function
$d^r\!\left(s,\widehat{s}\right)$ ($r\ge1$) replacing $d^2\!\left(s,\widehat{s}\right)$. In
such a case, the results we prove below for the quadratic risk easily extend to the risk 
$\Bbb{E}_s\cro{d^r(s,\widehat{s})}$. For simplicity, we shall only focus on the case $r=2$.
%

\subsection{Some illustrations\labs{F2}}
The previous theorem actually holds for various statistical frameworks. Let us provide a partial
list. 

\paragraph{Gaussian frameworks}
A prototype for Gaussian frameworks is provided by some Gaussian isonormal linear 
process as described in Section~2 of Birg\'e and Massart~\citeyearpar{MR1848946}. In such a 
case, $\bm{X}$ is a Gaussian linear process with a known  variance $\tau$, indexed by a subset
${\cal S}$ of some Hilbert space $\Bbb{L}_2(E,\mu)$. This means that $s\in{\cal S}$ determines
the distribution $P_s$. Regression with Gaussian errors and Gaussian sequences can both be seen
as particular cases of this framework. Then Theorem~\ref{T-models} is a consequence of
Theorem~2 of Birg\'e and Massart~\citeyearpar{MR1848946}.   In the regression setting, the
practical case of an unknown variance has been considered in Baraud, Giraud and
Huet~\citeyearpar{BaGiHu2009} under the additional assumption that
$\D(S)\vee\Delta_{\pi}(S)\le n/2$ for all $S\in\SS$. 

\paragraph{Density estimation}
Here $\bm{X}=(X_1,\ldots,X_n)$ is an $n$-sample with density $s^2$ with respect to
$\mu$ and ${\cal S}$ is the set of nonnegative elements of norm 1 in $\Bbb{L}_2(E,\mu)$.
Then $d(s,t)=\sqrt{2}h\left(s^2,t^2\right)$ where $h$ denotes the Hellinger distance
between densities,  $\tau=n^{-1}$ and Theorem~\ref{T-models} follows from Theorem~6
of Birg\'e~\citeyearpar{MR2219712}. Alternatively, one can take for $s$ the density itself, for 
${\cal S}$ the set of nonnegative elements of norm 1 in $\Bbb{L}_1(E,\mu)$ and set $q=1$. The
result then follows from Theorem~8 of Birg\'e~\citeyearpar{MR2219712}. Under the additional
assumption that $s\in\IL_{2}(E,\mu)\cap \IL_{\infty}(E,\mu)$, the case $q=2$ follows from
Theorem~6 of Birg\'e ~\citeyearpar{2008arXiv0808.1416B} with
$\tau=n^{-1}\norm{s}_{\infty}(1\vee\log\norm{s}_{\infty})$. 

\paragraph{Regression with fixed design}
We observe $\bm{X}=\{(x_1,Y_1),\ldots,(x_n,Y_n)\}$ with $\Bbb{E}[Y_i]=s(x_i)$ where $s$ is 
a function from $E=\{x_1,\ldots,x_n\}$ to $\R$ and the errors $\eps_{i}=Y_{i}-s(x_{i})$ are i.i.d.. Here $\mu$ is the uniform distribution on $E$, hence
$d^2(s,t)=n^{-1}\sum_{i=1}^n[s(x_i)-t(x_i)]^2$ and $\tau=1/n$. When the errors
$\eps_{i}$ are subgaussian, Theorem~\ref{T-models} follows from Theorem~3.1 in 
Baraud, Comte and Viennet~\citeyearpar{MR1845321}. For more heavy-tailed distributions (Laplace, Cauchy, etc.) we refer to Theorem~6 of Baraud~\citeyearpar{SMH} when $s$ takes its values in $[-1,1]$. 

\paragraph{Bounded regression with random design}
Let $(X,Y)$ be a pair of random variables with values in $E\times [-1,1]$ where $X$ has 
distribution $\mu$ and $\Bbb{E}[Y|X=x]=s(x)$ is a function from $E$ to $[-1,1]$. Our aim here
is to estimate $s$ from the observation of  $n$ independent copies 
$\bm{X}=\{(X_1,Y_1),\ldots,(X_n,Y_n)\}$ of $(X,Y)$. Here the distance $d$ corresponds to the
$\mathbb{L}_{2}(E,\mu)$-distance and Theorem~\ref{T-models} follows from Corollary~8 in
Birg\'e~\citeyearpar{MR2219712} with $\tau=n^{-1}$. 

\paragraph{Poisson processes}
In this case, $\bm{X}$ is a Poisson process on $E$ with mean measure    
$s^2\cdot\mu$, where $s$ is a nonnegative element of $\Bbb{L}_2(E,\mu)$. Then
$\tau=1$ and Theorem~\ref{T-models} follows from Birg\'e~\citeyearpar{Birge-Poisson}.
%

\Section{Approximation of  functions\labs{A}}
In this section, we give a brief overview of more or less classical collections of models commonly used for approximating smooth (and less smooth) functions on $[-1,1]^{k}$.  We shall
use these collections to approximate composite functions of the form  $g\circ  u$ by
approximating $g$ and $u$ separately. Finally, we shall explain the basic ideas of our approach
at the end of this section.

Collections of models with the property below will be of special interest throughout this paper. 
\begin{assumption}\lab{B-complexity}
The elements of the collection $\SS$ are finite dimensional linear spaces and for each $D\in
 \mathbb{N}$ the number of elements of $\SS$ with dimension $D$ is bounded by
$\exp[c(\SS)(D+1)]$ for some nonnegative constant $c(\SS)$ depending on $\SS$ only.
\end{assumption}

\subsection{Classical models for approximating functions}\label{sect-ApproxS}
Along this section, $d$ denotes the $\IL_{2}$-distance in $\IL_{2}([-1,1]^{k},2^{-k}dx)$, thus 
taking  $q=2$, $E=[-1,1]^{k}$ and $\mu$ the Lebesgue probability on $E$. 

\subsubsection{Approximating smooth functions on $[-1,1]^{k}$}
When $k=1$, a typical smoothness condition for a function $s$ on $[-1,1]$ is that it  belongs to 
some H\"older space ${\cal H}^{\alpha}([-1,1])$ with
$\alpha=r+\alpha'$, $r\in\Bbb{N}$ and $0<\alpha'\le1$ which is the set of all functions $f$ 
on $[-1,1]$ with a continuous derivative of order $r$ satisfying, for some constant $L(f)>0$,
\[
\left|f^{(r)}(x)-f^{(r)}(y)\right|\le L(f)|x-y|^{\alpha'}\quad\mbox{for all }x,y\in[-1,1].
\]
This notion of smoothness extends to functions $f(x_1,\ldots,x_k)$ defined on $[-1,1]^k$,
by saying that $f$ belongs to ${\cal H}^{\galpha}([-1,1]^{k})$ with 
$\bm{\alpha}=(\alpha_1,\ldots,\alpha_k)\in(0,+\infty)^k$ if, viewed as a function of $x_i$ 
only,  it belongs to ${\cal H}^{\alpha_{i}}([-1,1])$ for $1\le i\le k$ with some constant $L(f)$
independent of both $i$ and the variables $x_j$ for $j\ne i$. The smoothness of a function $s$
in ${\cal H}^{\galpha}([-1,1]^{k})$ is said to be {\it isotropic} if the $\alpha_{i}$ are all equal 
and {\it anisotropic} otherwise, in which case the quantity  $\overline{\alpha}$ given by
\[
\frac{1}{\overline{\alpha}}=\frac{1}{k}\sum_{i=1}^k\frac{1}{\alpha_i}
\] 
corresponds to the average smoothness of $s$. It follows from results in Approximation Theory that functions in the H\"older space 
${\cal H}^{\galpha}([-1,1]^{k})$ can be well approximated by piecewise polynomials on 
$k$-dimensional hyperrectangles. More precisely, our next proposition follows from results in 
Dahmen, DeVore and Scherer~\citeyearpar{MR581486}.
\begin{proposition}\lab{P-Approx1}
Let $(k,r)\in\N^{\star}\times\N$. There exists a collection of models $\HH_{k,r}=
\bigcup_{D\ge 1}\HH_{k,r}(D)$ satisfying Assumption~\ref{B-complexity} such that for each
positive integer $D$, the family  $\HH_{k,r}(D)$ consists of linear spaces $S$ with dimensions
${\cal D}(S)\le C'_1(k,r)D$ spanned by piecewise polynomials of degree at most $r$ on
$k$-dimensional hyperrectangles and for which 
\begin{equation}
\inf_{S\in \HH_{k,r}(D)}d(s,S)\le \inf_{S\in \HH_{k,r}(D)}d_{\infty}(s,S)\le
C'_2(k,r)L(s)D^{-\overline{\alpha}/k},
\labe{Eq-smooth}
\end{equation}
for all $s\in {\cal H}^{\galpha}([-1,1]^{k})$ with $\sup_{1\le i\le k}\alpha_i\le r+1$.
\end{proposition}

\subsubsection{Approximating functions in anisotropic Besov spaces}
Anisotropic Besov spaces generalize anisotropic H\"older spaces and are defined in a similar
way by using directional moduli of smoothness, just as  H\"older spaces are defined using
directional derivatives. To be short, a function belongs to an anisotropic Besov space on
$[-1,1]^k$ if, when all coordinates are fixed apart from one, it belongs to a Besov space on
$[-1,1]$. A precise definition (restricted to $k=2$ but which can be generalized easily) can be
found in Hochmuth~\citeyearpar{MR1884234}. The general definition together with useful
approximation properties by piecewise polynomials can be found in
Akakpo~\citeyearpar{Akakpo09}. For $0<p\le+\infty$, $k>1$ and $\bm{\beta}\in
(0,+\infty)^{k}$, let us denote by $\B_{p,p}^{\bm{\beta}}([-1,1]^k)$ the  anisotropic Besov
spaces. In particular, $\B_{\infty,\infty}^{\bm{\beta}}([-1,1]^k) ={\cal H}^{\gbeta}([-1,1]^{k})$.
It follows from Akakpo~\citeyearpar{Akakpo09} that Proposition~\ref{P-Approx1} can be
generalized to Besov spaces in the following way.
%
\begin{proposition}\lab{P-Approx2}
Let $p>0$, $k\in \N^{\star}$ and $r\in\N$. There exists a collection of models 
$\BB_{k,r}=\bigcup_{D\ge 1}\BB_{k,r}(D)$ satisfying Assumption~\ref{B-complexity} such
that for each positive integer $D$, $\BB_{k,r}(D)$ consists of linear spaces $S$ with dimensions
${\cal D}(S)\le C'_1(k,r)D$ spanned by piecewise polynomials of degree at most $r$ on
$k$-dimensional hyperrectangles and for which 
\begin{equation}
\inf_{S\in \BB_{k,r}(D)}d(s,S)\le
C'_2\!\left(k,r,p\right)|s|_{\bm{\beta},p,p}D^{-\overline{\beta}/k}
\labe{Eq-smooth}
\end{equation}
for all $s\in {\cal B}_{p,p}^{\gbeta}([-1,1]^{k})$ with semi-norm  $|s|_{\bm{\beta},p,p}$ and
$\gbeta$ satisfying
\begin{equation}
\sup_{1\le i\le k}\beta_i<r+1\qquad\mbox{and}\qquad
\overline{\beta}>k\left[\left(p^{-1}-2^{-1}\right)\vee0\right].
\labe{Eq-Bes1}
\end{equation}
\end{proposition}
%

\subsection{Approximating composite functions\labs{A2}}

\subsubsection{Preliminary approximation results\labs{A2a}}
As already mentioned, we assume along the paper that $s$ is either of the form $g\circ u$ or can 
be well approximated by a function of this form, where $u=(u_{1},\ldots,u_{l})$ belongs to
$\T^{l}$ and $g$ to $\F_{l,c}$. The purpose of this section is to see how well $f\circ t$
approximates $g\circ u$ when we know how well $f$ approximates $g$ and
$t=(t_{1},\ldots,t_{l})$ approximates $u$. We start with the definition of the modulus of 
continuity of a function $g$ in $\F_{l,c}$.
%
\begin{definition}\lab{D-modcont}
We say that $\w$ from $[0,2]^{l}$ to $\R_{+}^{l}$ is a modulus of continuity for a continuous 
function $g$ on $[-1,1]^{l}$ if for all $z\in [0,2]^{l}$, $\w(z)$ is of the form
$\w(z)=(\w_1(z_{1}),\ldots,\w_l(z_{l}))$ where each function $\w_{j}$ with $j=1,\ldots,l$ is
continuous, nondecreasing and concave from $[0,2]$ to $\Bbb{R}_+$, satisfies $\w_j(0)=0$, and 
\[
|g(x)-g(y)|\le\sum_{j=1}^l \w_j(|x_j-y_j|)\quad\mbox{for all }x,y\in[-1,1]^l.
\]
For $\galpha\in(0,1]^{l}$ and $\gL\in(0,+\infty)^l$, we say that $g$ is 
$(\galpha,\gL)$-H\"olderian if one can take $\w_{j}(z)=L_{j}z^{\alpha_{j}}$ for all $z\in [0,2]$
and $j=1,\ldots,l$. It is said to be $\gL$-Lipschitz if it is $(\galpha,\gL)$-H\"olderian with
$\galpha=(1,\ldots,1)$.
\end{definition}
Note that our definition of a modulus of continuity implies that the $\w_{j}$ are subadditive, a 
property which we shall often use in  the sequel and that, given $g$, one can always choose for
$\w_j$ the least concave majorant of $w_j$ where
\[
w_j(z)=\sup_{x\in[-1,1]^l;\; x_j\le1-z}
|g(x)-g(x_1,\ldots,x_{j-1},x_j+z,x_{j+1},\ldots,x_l)|.
\]
Then $\w_j(z)\le2w_j(z)$, according to Lemma~6.1 p.~43 of DeVore and
Lorentz~\citeyearpar{DeVore}.

Let us now turn to our main approximation result to be proved in Section~\ref{H6}.
%
%
\begin{proposition}\lab{P-Approx}
Let $p\ge 1$, $g\in\F_{l,c}$, $f\in\F_{l,\infty}$ and $t,u\in\T^{l}$. If $\w_g$ is a modulus of continuity for $g$, 
then
\[
\|g\circ u-f\circ t\|_{p}\le
d_{\infty}(g,f)+2^{1/p}\sum_{j=1}^{l}\w_{g,j}\left(\|u_{j}-t_{j}\|_{p}\right)
\]
with the convention $2^{1/\infty}=1$.
\end{proposition}
%

\subsubsection{The main ideas underlying our construction\labs{A2b}}
Let us take $p=q=2$ and $E=[-1,1]^{k}$ with $k>l\ge 1$. A consequence of 
Proposition~\ref{P-Approx} is the following. If one considers a finite dimensional linear space
$F\subset \F_{l,\infty}$ for approximating $g$ and compact sets $T_{j}\subset \T$ for
approximating the $u_{j}$, there exists $t\in \gT= \prod_{j=1}^{l}T_{j}$ such that the linear
space $S_{t}=\{f\circ t\,|\,f\in F\}$ approximates the composite function $g\circ u$  with an
error bound 
\begin{equation}\label{JB007}
d(g\circ u,S_{t})\le d_{\infty}(g,F)+\sqrt{2}\sum_{j=1}^{l}\w_{g,j}\left(d(u_{j},T_{j})\right).
\end{equation}
The case where the function $g$ is Lipschitz is of particular interest since, up to constants,  the
error bound we get is the sum of those for approximating separately $g$ by $F$ (with respect to
the $\IL_\infty$-distance) and the
$u_{j}$ by $T_{j}$. In particular, if $s$ were exactly of the form $s=g\circ u$ for some known
functions
$u_{j}$, we could use a suitable linear space $F$ with dimension of order $D$ in $\HH_{l,0}(D)$ to
approximate $g$, and take $T_{j}=\{u_{j}\}$ for all $j$. In this case the linear space $S_{u}$
whose dimension is also of order $D$ would approximate $s=g\circ u$ with an error bounded
by $D^{-1/l}$. Note that if the $u_{j}$ did belong to some H\"older space
$\H^{\gbeta}([-1,1]^{k})$ with $\gbeta\in(0,1]^{k}$, the overall regularity of the function
$s=g\circ u$ could not be expected to be better than $\gbeta$-H\"olderian, since this regularity
is already achieved by taking $g(y_{1},\ldots,y_{l})=y_{1}$. In comparison, an approach
based on the overall smoothness of $s$, which would completely ignore the fact that $s=g\circ
u$ and the knowledge of the $u_{j}$, would lead to an approximation bound of order
$D^{-\overline \beta/k}$. The former bound, $D^{-1/l}$, based on the structural assumption
that $s=g\circ u$ therefore improves on the latter since $\overline \beta\le 1$ and $k>l$. Of
course, one could argue that the former approach uses the knowledge of the $u_{j}$, which is
quite a strong assumption for statistical issues. Actually, a more reasonable approach would be
to assume that $u$ is unknown but close to a parametric set $\overline\gT$, in which case, it
would be natural to replace the single model $S_{u}$ used for approximating $s$, by the
family of models $\SS_{\overline \gT}(F)=\{S_{t}\,|\,t\in\overline\gT\}$  and, ideally,  let the
usual model selection techniques select some best linear space among it. Unfortunately, results
such as Theorem~\ref{T-models} do not apply to this case, since the family $\SS_{\overline
\gT}(F)$ has the same cardinality as $\overline \gT$ and is therefore typically not countable.
The main idea of our approach is to take advantage of the fact that the $u_{j}$ take their values in
$[-1,1]$ so that we can embed $\overline \gT$ into a compact subset of $\T^{l}$. We may then
introduce a suitably discretized version $\gT$ of $\overline \gT$ (more precisely, of its
embedding) and replace the ideal collection $\SS_{\overline \gT}(F)$ by $\SS_{\gT}(F)$, for
which similar approximation properties can be proved. The details of this discretization device
will be  given in the proofs of our main results. Finally, we shall let both $\overline \gT$ and
$F$ vary into some collections of models and use all the models of the various resulting
collections $\SS_{\gT}(F)$ together in order to estimate $s$ at best.   

\Section{The basic theorems\labs{B}}

\subsection{Model selection using classical approximation spaces\labs{B1}}
If we assume that the unknown parameter $s$ to be estimated is equal or close to some 
composite function of the form $g\circ u$ with $u\in\T^{l}$ and $g\in\F_{l,c}$ and if we wish 
to estimate $g\circ u$ by model selection we need to have at disposal families of models for
approximating both $g$ and the components $u_j$, $1\le j\le l$ of $u$. As already seen in 
Section~\ref{sect-ApproxS}, typical sets that are  used for approximating elements of $\F_{l,c}$
or $\T^l$ are finite-dimensional linear spaces or subsets of them and we need a theorem
which applies to such classical approximation sets for which it will be convenient to choose
the following definition of dimension.
%
\begin{definition}\lab{D-dimension}
Let $H$ be a linear space and $S\subset H$. The dimension $\D(S)\in\N\cup\{\infty\}$ of 
$S$ is 0 if $|S|=1$ and is, otherwise, the dimension (in the usual sense) of the linear span of 
$S$.
\end{definition}
Our construction of  estimators $\widehat s$ of $g\circ u$ will be based on some set $\frak{S}$
of the following form:
\begin{equation}
\frak{S}=\{l,\FF,\gamma,\TT_1,\ldots,\TT_l,\lambda_1,\ldots,\lambda_l\},\quad
l\in\Bbb{N}^\star,
\labe{Eq-frakS}
\end{equation}
where $\FF,\TT_1,\ldots,\TT_l$ denote families of models and $\gamma$, $\lambda_j$ are
measures on $\FF$ and $\TT_j$ respectively. In the sequel, we shall assume that $\frak{S}$
satisfies the following requirements.
%
\begin{assumption}\lab{A-models}
The set $\frak{S}$ is such that

i) the family $\FF$ is a countable set and consists of finite-dimensional linear subspaces $F$
of  $\F_{l,\infty}$ with respective dimensions $\D(F)\ge1$,

ii) for $j=1,\ldots,l$, $\TT_{j}$ is a countable set of subsets of  $\mathbb{L}_q(E,\mu)$ with
finite dimensions,

iii) the measure $\gamma$ is a subprobability on $\FF$,

iv) for $j=1,\ldots,l$, $\lambda_j$ is a subprobability on $\TT_j$.
\end{assumption}
Given $\frak{S}$, one can design an estimator $\widehat s$ with the following properties.
%
\begin{theorem}\lab{T-general}
Let $\frak{S}$ satisfy Assumption~\ref{A-models}. One can build an estimator
$\widehat{s}=\widehat s(\bm{X})$ satisfying, for all $u\in\T^l$ and $g\in\F_{l,c}$ with
modulus of continuity $\w_g$,
\begin{eqnarray}
C\Bbb{E}_s\!\left[d^2\left(s,\widehat s)\right)\st\right]&\le&\sum_{j=1}^l \inf_{T\in
\TT_j}\ac{l\w^2_{g,j}\!\left(d(u_j,T)\st\right)+\tau\cro{\Delta_{\lambda_{j}}(T)+i(g,j,T)
\D(T)}}\nonumber\\&&\mbox{}+ d^2(s,g\circ u)+\inf_{F\in\FF}\ac{d_{\infty}^{2}(g,F)+
\tau\cro{\Delta_{\gamma}(F)+\D(F)}}\!,\qquad
\labe{Eq-MainRB}
\end{eqnarray}
where $i(g,j,T)=1$ if ${\cal D}(T)=0$ and otherwise,
\begin{equation}
i(g,j,T)=\inf\left\{i\in\N^\star\,|\,l\w^2_{g,j}\left(e^{-i}\right)\le\tau i{\cal D}(T)\right\}
<+\infty.
\labe{Eq-igj}
\end{equation}
\end{theorem}
Note that, since the risk bound~(\ref{Eq-MainRB}) is valid for all $g\in\F_{l,c}$ and
$u\in\T^l$, we can minimize the right-hand side of~(\ref{Eq-MainRB}) with respect to $g$ and
$u$ in order to optimize the bound. The proof of this theorem is postponed to Section~\ref{H3}.

Of special interest is the case where $g$ is $\gL$-Lipschitz. If one is mainly interested by the
dependence of the risk bound with respect  to $\tau$ as it tends to 0, one can check that
$i(g,j,T)\le\log\tau^{-1}$ for $\tau$ small enough (depending on $l$ and $\gL$) so that
(\ref{Eq-MainRB}) becomes for such a small $\tau$
\begin{eqnarray*}
C'\Bbb{E}_s\!\left[d^2\left(s,\widehat s\right)\st\right]&\le&\sum_{j=1}^l\inf_{T\in \TT_j}
\ac{d^{2}(u_j,T)+\tau\pa{\Delta_{\lambda_j}(T)+\D(T)\log\tau^{-1}}}\\&&+ d^2(s,g\circ
u)+\inf_{F\in\FF}\ac{d_{\infty}^{2}(g,F)+\tau\cro{\D(F)+\Delta_{\gamma}(F)}}.
\end{eqnarray*}
If it were possible to apply Theorem~\ref{T-models} to the models $F$ with the distance
$d_\infty$ and the models $T$ with the distance $d$ for each $j$ separately, we would get risk 
bounds of this form, apart from the value of $C'$ and the extra $\log\tau^{-1}$ factor. This
means that, apart from this extra logarithmic factor, our procedure amounts to performing $l+1$
separate model selection procedures, one with the collection $\FF$ for estimating $g$ and the
other ones with the collections $ \TT_{j}$ for the components $u_{j}$ and finally getting the sum
of  the $l+1$ resulting risk bounds. The result is however slightly different when $g$ is no 
longer Lipschitz. When $g$ is (\galpha,\gL)-H\"olderian then one can check that $i(g,j,T)\le
\L_{j,T}$ where $\L_{j,T}=1$ if $\D(T)=0$ and, if $\D(T)\ge1$,
\begin{eqnarray}
\L_{j,T}&=&\left[\alpha_j^{-1}\log\left(lL_j^2[\tau\D(T)]^{-1}\right)\right]\bigvee1
\labe{Eq-LJT}\\ 
&\le& C'(l,\alpha_{j})\cro{\log(\tau^{-1})\vee \log(L_{j}^{2}/\D(T))\vee1}.
\labe{Eq-LJT'}
\end{eqnarray}
In this case, Theorem~\ref{T-general} leads to the following result.
\begin{corollary}\label{C-Main0}
Assume that the assumptions of Theorem~\ref{T-general} holds. For all $(\galpha,\gL)$-H\"olderian function $g$ with $\galpha\in(0,1]^{l}$ and $\gL\in(\R_{+}^{\star})^{l}$, the estimator $\widehat s$ of Theorem~\ref{T-general} satisfies
\begin{eqnarray}
C\Bbb{E}_s\!\left[d^2\left(s,\widehat s\right)\st\right]&\le&\nonumber
\sum_{j=1}^l\inf_{T\in\TT_{j}}\ac{lL_{j}^{2}d^{2\alpha_{j}}(u_j,T)+
\tau\left[\Delta_{\lambda_{j}}(T)+\D(T)\L_{j,T}\right]}\\&&\mbox{}+d^2(s,g\circ u)+
\inf_{F\in\FF}\ac{d_{\infty}^{2}(g,F)+\tau\cro{\Delta_{\gamma}(F)+\D(F)}}\!,\qquad
\labe{Eq-RB2a}
\end{eqnarray}
where $\L_{j,T}$ is defined by~\eref{Eq-LJT} and bounded by ~\eref{Eq-LJT'}
\end{corollary}
%
\subsection{Mixing collections corresponding to different values of $l$\labs{B4}}
If it is known that $s$ takes the special form  $g\circ u$ for some unknown values of
$g\in\F_{l,c}$ and $u\in \T^l$, or if $s$ is very close to some function of the form $g\circ
u$, the previous approach is quite satisfactory. If we do not have such an information, we may
apply the previous construction with several values of $l$ simultaneously, approximating $s$ 
by different combinations $g_l\circ u_l$ with $u_l$ taking its values in $[-1,1]^l$, $g_l$ a function
on $[-1,1]^l$ and $l$ varying among some subset $I$ of $\N^{\star}$. To each value of $l$, we
associate, as before,  $l+1$ collections of models and the corresponding subprobabilities, each
$l$ then leading to an estimator $\widehat{s}_l$ the risk of which is bounded by
$\risq(\widehat s_{l},g_{l},u_{l})$ given by the right-hand side of (\ref{Eq-MainRB}). The model
selection approach allows us to use all the previous collections of models for all values of $l$
simultaneously to build a new estimator the risk of which is approximately as good as the risk
of the best of the $\widehat{s}_l$. More generally, let us assume that we have at hand a countable
family $\{\frak{S}_{\ell},\,\ell\in I\}$ of sets $\frak{S}_{\ell}$ of the form (\ref{Eq-frakS})
satisfying Assumption~\ref{A-models} for some $l=l({\ell})\ge 1$. To each such set,
Theorem~\ref{T-general} associates an estimator $\widehat s_{\ell}$ with a risk bounded by
\[
\Bbb{E}_s\!\left[d^2\left(s,\widehat{s}_{\ell}\right)\st\right]\le\inf_{(g,u)}\risq(\widehat
s_{\ell}, g,u),
\]
where $\risq(\widehat s_{\ell},g,u)$ denotes the right-hand side of (\ref{Eq-MainRB}) when
$\frak{S}=\frak{S}_{\ell}$ and the infimum runs among all pairs $(g,u)$ with $g\in\F_{l(\ell),
c}$ and $u\in\T^{l(\ell)}$. We can then prove (in Section~\ref{H4} below) the following result.
%
\begin{theorem}\lab{T-melange}
Let $I$ be a countable set and $\nu$ a subprobability on $I$. For each $\ell\in I$ we are
given a set $\frak{S}_{\ell}$ of the form (\ref{Eq-frakS}) that satisfies
Assumption~\ref{A-models} with $l=l({\ell})$ and a corresponding estimator $\widehat{s}_l$
provided by Theorem~\ref{T-general}. One can then design a new estimator
$\widehat{s}=\widehat{s}(\bm{X})$ satisfying 
\[
C\Bbb{E}_s\!\left[d^2\left(s,\widehat{s}\right)\st\right]\le \inf_{\ell\in I}\inf_{(g,u)}
\ac{\risq(\widehat s_{\ell},g,u)+\tau\Delta_{\nu}(\ell)},
\]
where $\risq(\widehat s_{\ell},g,u)$ denotes the right-hand side of (\ref{Eq-MainRB}) when
$\frak{S}=\frak{S}_{\ell}$ and the second infimum runs among all pairs $(g,u)$ with
$g\in\F_{l(\ell),c}$ and $u\in\T^{l(\ell)}$.
\end{theorem}

\Section{Applications\labs{C}}
The aim of this section is to provide various applications of Theorem~\ref{T-general} and its 
corollaries.  For approximating functions on $[-1,1]^{k}$, we shall repeatedly use  the families
$\HH_{k,r}$ and $\BB_{k,r}$ introduced in Section~\ref{sect-ApproxS}. Since for all $r\ge
0$, $\HH_{k,r}$ satisfies Assumption~\ref{B-complexity} for some constant $c(\HH_{k,r})$,
the measure $\gamma$ on $\HH_{k,r}$ defined by 
\begin{equation}
\Delta_{\gamma}(S)=(c(\HH_{k,r})+1)(D+1)\quad\mbox{for all }S\in\HH_{k,r}(D)\setminus
\bigcup_{1\le D'<D}\HH_{k,r}(D')
\labe{Eq-defgamma}
\end{equation}
is a subprobability since
\[
\sum_{S\in\HH_{k,r}}e^{-\Delta_{\gamma}(S)}\le\sum_{D\ge 1}e^{-D}
\ab{\st\HH_{k,r}(D)}e^{-c(\HH_{k,r})(D+1)}\le \sum_{D\ge 1}e^{-D}<1.
\]
We shall similarly consider the subprobability $\lambda$ defined on $\BB_{k,r}$ by
\begin{equation}
\Delta_{\lambda}(S)=(c(\BB_{k,r})+1)(D+1)\quad\mbox{for all }S\in\BB_{k,r}(D)
\setminus\bigcup_{1\le D'<D}\BB_{k,r}(D').
\labe{Eq-deflambda}
\end{equation}
Finally, for $g\in \H^{\galpha}([-1,1]^{l})=\B^{\galpha}_{\infty,\infty}([-1,1]^{l})$ with 
$\galpha\in(\R_{+}^{\star})^{l}$,  we set
$\norm{g}_{\galpha,\infty}=|g|_{\galpha,\infty,\infty}+\inf L'$ where the infimum runs
among all numbers $L'$ for which $\w_{g,j}(z)\le L'z^{\alpha_{j}\wedge 1}$ for all $z\in
[0,2]$ and $j=1,\ldots,l$.

\subsection{Estimation of smooth functions on $[-1,1]^k$\labs{C1}}
In this section, our aim is to establish risk bounds for our estimator $\widehat s$ when 
$s=g\circ u$ for some smooth functions $g$ and $u$. We shall discuss the improvement, in
terms of rates of convergence as $\tau$ tends to 0, when assuming such a structural
hypothesis, as  compared to a pure smoothness assumption on $s$. Throughout this section,
we take $q=2$,
$E=[-1,1]^{k}$ and $d$ as the $\IL_{2}$-distance on $\IL_{2}(E,2^{-k}dx)$.

\subsubsection{Convergence rates using composite functions\labs{C1c}}
Let us consider here the set $\S_{k,l}(\galpha,\gbeta,\gp,L,\gR)$ gathering the composite 
functions $g\circ u$ with $g\in\H^{\galpha}([-1,1]^{l})$ satisfying
$\norm{g}_{\galpha,\infty}\le L$ and $u_{j}\in \B_{p_{j},p_j}^{\gbeta_{j}}$ with
semi-norms $\ab{u_{j}}_{\gbeta_{j},p_{j},p_j}\le R_{j}$ for all $j=1,\ldots,l$. The following
result holds. 
%
\begin{theorem}\lab{T-Smooth}
There exists an estimator $\widehat s$ such that, for all $l\ge 1$, $\galpha,\gR \in 
(\R_{+}^{\star})^{l}$, $L>0$, $\gbeta_{1},\ldots,\gbeta_{l}\in(\R_{+}^{\star})^{k}$ and
$\gp\in(0,+\infty]^{l}$ with 
$\overline{\beta}_j> k\left[\left(p_j^{-1}-2^{-1}\right)\vee0\right]$ for $1\le j\le l$,  
\begin{eqnarray*}
\lefteqn{\sup_{s\in\S_{k,l}(\galpha,\gbeta,\gp,L,\gR)}C'\Bbb{E}_s\!\left[d^2\!
\left(s,\widehat{s}\right)\st\right]}\hspace{20mm}\\ &\le&\sum_{j=1}^{l}
\left(LR_{j}^{\alpha_{j}\wedge 1}\right)^{{2k\over 2\overline{\beta_j}(\alpha_j\wedge1)+
k}}\left[\tau\L\right]^ {{2\overline{\beta_j}(\alpha_j\wedge1)\over
2\overline{\beta_j}(\alpha_j\wedge1)+k}}\ +\   L^{{2l\over
l+2\overline{\alpha}}}\tau^{{2\overline{\alpha}\over l+2\overline{\alpha}}} +\tau\L,
\end{eqnarray*}
where $\L=\log(\tau^{-1})\vee\log(L^{2})\vee 1$ and $C'$ depends on $k,l,\bm{\galpha}$,
$\bm{\beta}$ and $\gp$.
\end{theorem}
Let us recall that we need not assume that $s$ is exactly of the form $g\circ u$ but rather, as 
we  did  before, that $s$ can be approximated by a function $\overline{s}=g\circ u\in
\S_{k,l}(\galpha,\gbeta,\gp,L,\gR)$. In such a case we simply get an additional bias term 
of the form $d^2\left(s,\overline{s}\right)$ in our risk bounds. \vspace{2mm}\\
%
\noindent{\em Proof:}
Let us fix some value of $l\ge 1$ and take $s=g\circ
u\in\S_{k,l}(\galpha,\gbeta,\gp,L,\gR)$ and define 
\[
r=r(\galpha,\gbeta)=1+\left\lfloor \max_{i=1,\ldots,l}\alpha_{i}\bigvee\max_{j=1,\ldots,l,
\ell=1,\ldots,k}\beta_{j,\ell}\right\rfloor.
\]
The regularity properties of $g$ and the $u_{j}$ together with Propositions~\ref{P-Approx1} 
and~\ref{P-Approx2} imply that for all $D\ge 1$, there exist $F\in\HH_{l,r}(D)$ and sets
$T_{j}\in\BB_{k,r}(D)$ for $j=1,\ldots,l$ such that 
\[
{\cal D}(F)\le C'_1\!\left(l,\bm{\alpha},\gbeta\right)D;\qquad
d_{\infty}(g,F)\le C_2'\!\left(l,\bm{\alpha},\gbeta\right)LD^{-\overline{\alpha}/l};
\]
and, for $1\le j\le l$,
\[
{\cal D}(T_j)\le C'_3\!\left(k,\galpha,\bm{\beta}_j,p_j\right)D;\qquad d(u_j,T_j)\le
C'_4\!\left(k,\galpha,\bm{\beta}_j,p_j\right)R_{j}D^{-\overline{\beta_j}/k}.
\]
Since the collections $\HH_{l,r}$ and $\BB_{k,r}$ satisfy Assumption~\ref{B-complexity} 
and $\w_{g,j}(z)\le Lz^{\alpha_j\wedge1}$ for all $j$ and $z\in[0,2]$, we may apply
Corollary~\ref{C-Main0} with 
\[
\frak{S}_{l,r}=(l,\HH_{l,r},\gamma_{r},\BB_{k,r},\ldots,\BB_{k,r},\lambda_{r},\ldots,
\lambda_{r})
\]
the subprobabilities $\gamma_{l,r}$ and $\lambda_{l,r}$ being given
by~\eref{Eq-defgamma}  and~\eref{Eq-deflambda} respectively. Besides, it follows from 
(\ref{Eq-LJT'}) that $\L_{j,T}\le C'(l,\galpha)\L$ for all $j$, so that \eref{Eq-RB2a} implies
that the risk of the resulting estimator $\widehat s_{l,r}$ is bounded from above by  
\[
C'\RR(\widehat s_{l,r},g,u)=\sum_{j=1}^l\inf_{D\ge 1}\cro{L^{2}R_{j}^{2(\alpha_{j}
\wedge 1)}D^{-2(\alpha_{j}\wedge 1)\overline \beta_{j}/k}+D\tau\L}+\inf_{D\ge
1}\cro{L^{2}D^{-2\overline{\alpha}/l}+D\tau},
\]
for some constant $C'$ depending on $l,k,\galpha,\gbeta_{1},\ldots,\gbeta_{l}$. We obtain 
the result by optimizing each term of the sum with respect to $D$, by means of
Lemma~\ref{L-opt}, and by using Theorem~\ref{T-melange} with $\nu$ defined for
$\ell=(l,r)\in\N^{\star}\times\N$ by $\nu(l,r)=e^{-(l+r+1)}$ for which
$\Delta_{\nu}(l,r)\tau\le (l+r+1)\RR(\widehat s_{l,r},g,u)$ for all $l,r$.\cqfd
%

\subsubsection{Structural assumption versus smoothness assumption\labs{C1a}}
In view of discussing the interest of the risk bounds provided by Theorem~\ref{T-Smooth}, let
us focus here, for simplicity, on the case where $g\in{\cal H}^{\alpha}([-1,1])$ with $\alpha>0$
(hence $l=1$) and $u$ is a function from $E=[-1,1]^k$ to $[-1,1]$ that belongs to 
${\cal H}^{\bm{\beta}}([-1,1]^{k})$ with $\gbeta\in(\R_{+}^{\star})^{k}$. The following
proposition is to be proved in Section~\ref{H7}. 
%
\begin{proposition}\lab{}
\lab{P-Holder}
Let $\phi$ be the function defined on $(\R_{+}^{\star})^{2}$ by 
\[
\phi(x,y)=\left\{\begin{array}{ll}xy&\mbox{if }x\vee y\le1;\\
x\wedge y&\mbox{otherwise.}\end{array}\right.
\]
For all $k\ge 1$, $\alpha>0$, $\bm{\beta}\in(\R_{+}^{\star})^{k}$, $g\in{\cal H}^{\alpha}([-1,1])$ and $u\in{\cal H}^{\gbeta}([-1,1]^{k})$,
\begin{equation}
g\circ u\in{\cal H}^{\gtheta}([-1,1]^{k})\quad\mbox{with}\quad\theta_i=
\phi(\beta_i,\alpha)\quad\mbox{for }1\le i\le k. 
\labe{Eq-compsmooth}
\end{equation}
Moreover, $\bm{\theta}$ is the largest possible value for which (\ref{Eq-compsmooth}) holds 
for all $g\in{\cal H}^{\alpha}([-1,1])$ and $u\in{\cal H}^{\gbeta}([-1,1]^{k})$ since, whatever
$\bm{\theta}'\in (\R_{+}^{\star})^{k}$ such that $\theta'_{i}>\theta_{i}$ for some 
$i\in\ac{1,\ldots,k}$, there exists some $g\in{\cal H}^{\alpha}([-1,1])$ and $u\in 
{\cal H}^{\gbeta}([-1,1]^{k})$ such that $g\circ u\not\in{\cal H}^{\gtheta'}([-1,1]^{k})$.
\end{proposition}
Using the information that $s$ belongs to ${\cal H}^{\gtheta}([-1,1]^{k})$ with $\gtheta$ given
by (\ref{Eq-compsmooth}) and that we cannot assume that $s$ belongs to some smoother class
(although this may happen in special cases) since $\bm{\theta}$ is minimal, but ignoring the
fact that $s=g\circ u$, we can estimate $s$ at rate $\tau^{2\overline{\theta}/
(2\overline{\theta}+k)}$ (as $\tau$ tends to 0)  while, on the other hand, by using
Theorem~\ref{T-Smooth} and the structural information that $s=g\circ u$, we can achieve the
rate 
\[
\tau^{2\alpha/\left(2\alpha+1\right)}+\left(\tau\left[\log\tau^{-1}\right]\st
\right) ^{2\overline{\beta}(\alpha\wedge1)/(2\overline{\beta}(\alpha\wedge1)+k)}.
\]
Let us now compare these two rates. First note that it follows from (\ref{Eq-compsmooth}) that 
$\theta_i\le\alpha$ for all $i$, hence $\overline{\theta}\le\alpha$ and, since $k>1$,
$2\alpha/(2\alpha+1)>2\overline{\theta}/\left(2\overline{\theta}+k\right)$. Therefore the
term $\tau^{2\alpha/(2\alpha+1)}$ always improves over $\tau^{2\overline{\theta}/
\left(2\overline{\theta}+k\right)}$ when $\tau$ is small and, to compare the two rates, it is
enough to compare $\overline{\theta}$ with $\overline{\beta}(\alpha\wedge1)$. To do so, we
use the following lemma (to be proved in Section~\ref{H8}). 
\begin{lemma}\lab{L-comphold}
For all $\alpha>0$ and $\bm{\beta}\in(\R_{+}^{\star})^{k}$, the smoothness index 
\[
\bm{\theta}=\left(\st\phi(\alpha,\beta_{1}),\ldots,\phi(\alpha,\beta_{k})\right)
\] 
satisfies $\overline{\theta}\le\overline{\beta}(\alpha\wedge1)$ and equality holds if and only
if $\sup_{1\le i\le k}\beta_i\le\alpha\vee1$.
\end{lemma}
When $\sup_{1\le i\le k}\beta_i\le\alpha\vee1$, our special strategy does not bring any
improvement as compared to the standard one, it even slightly deteriorates the  risk bound
because of the extra $\log\tau^{-1}$ factor. On the opposite, if $\sup_{1\le i\le k}\beta_i>
\alpha\vee1$, our new strategy improves over the classical one and this improvement can be
substantial if $\overline{\beta}$ is much larger than $\alpha\vee1$. If, for instance, 
$\alpha=1$ and $\overline{\beta}=k=\beta_{j}$ for all $j$, we get a bound of order
$\left[\st\tau\left(\log\tau^{-1}\right)\right]^{2/3}$ which, apart from the extra
$\log\tau^{-1}$ factor, corresponds to the  minimax rate of estimation of a Lipschitz
function on $[-1,1]$, instead of the risk bound $\tau^{2/(2+k)}$ that we would get if we
estimated $s$ as a Lipschitz function on $[-1,1]^k$. When our strategy does not improve over 
the classical one, i.e.\ when $\sup_{1\le i\le k}\beta_i\le\alpha\vee1$, the additional loss due
to the extra logarithmic factor in our risk bound can be avoided by mixing the models used for
the classical strategy with the models used for designing our estimator, following the recipe of
Section~\ref{B4}.

\subsection{Generalized additive models\labs{C4}}
In this section, we assume that $E=[-1,1]^{k}$, $\mu$ is the Lebesgue probability on $E$ and 
$q=2$. A special structure that has often been considered in regression corresponds to functions
$s=g\circ u$ with
\begin{equation}\label{formadd}
u(x_1,\ldots,x_{k})=u_{1}(x_{1})+\ldots+u_{k}(x_{k});\quad
s(x)=g\left(\st u_1(x_1)+\ldots+u_k(x_k)\right),
\end{equation}
where the $u_{j}$ take their values in $[-1/k,1/k]$ for all $j=1,\ldots,k$. Such a model has been 
considered in Horowitz and Mammen~\citeyearpar{MR2382659} and while their approach is
non-adaptive, ours, based  on Theorem~\ref{T-general} and a suitable choice of the collections of
models, allows to derive a fully adaptive estimator with respect to the regularities of  $g$ and the 
$u_{j}$. More precisely, for $r\in\N$, let $\TT_{r}$ be the collection of all models of the form
$T=T_{1}+\ldots+T_{k}$ where for $j=1,\ldots,k$, $T_{j}$ is the set of functions of the form
$x\mapsto t_{j}(x_{j})$ with $x\in E$ and $t_{j}$ in $ \BB_{1,r}$. Using $\lambda_{r}=\lambda$
as defined by~\eref{Eq-deflambda}, we endow $\TT_{r}$ with the subprobability
$\lambda_{r}^{(k)}$ defined for $T\in \TT_{r}$ by the infimum of the quantities
$\prod_{i=1}^{k}\lambda_{r}(T_{i})$ when $(T_{1},\ldots,T_{k})$ runs among all the $k$-uplets
of $\BB_{1,r}^{k}$ satisfying $T=T_{1}+\ldots+T_{k}$. Finally, for $\alpha,L>0$,
$\gbeta,\gR\in(\R_{+}^{\star})^{k}$ and $\gp=(p_{1},\ldots,p_{k})\in(\R_{+}^{\star})^{k}$,
let $\S_{k}^{\rm Add}(\alpha,\gbeta,\gp,L,\gR)$ be the set of  functions of the
form~\eref{formadd} with $g\in\H^{\alpha}([-1,1])$ satisfying $\norm{g}_{\alpha,\infty}\le
L$ and $u_{j}\in\B^{\beta_{j}}_{p_{j},p_j}([-1,1])$ with $\ab{u_{j}}_{\beta_{j},p_{j},p_j}\le
R_{j}k^{-1}$ for all $j=1,\ldots,k$. Using the sets
$\frak{S}_{r}=(1,\HH_{1,r},\gamma_{r},\TT_{r},\lambda_{r}^{(k)})$ with $r\in\N$ we can
build an estimator with the following property. 
\begin{theorem}\label{additif}
There exists an estimator $\widehat s$ which satisfies for all $\alpha,L>0$, 
$\gp,\gR\in(\R_{+}^{\star})^{k}$ and $\gbeta\in(\R_{+}^{\star})^{k}$ with $\beta_{j}> 
\pa{1/p_{j}-1/2}_{+}$ for all $j=1,\ldots,k$,
\begin{eqnarray*}
\lefteqn{\sup_{s\in \S_{k}^{\rm Add}(\alpha,\gbeta,\gp,L,\gR)}C'\E_{s}
\cro{d^{2}(s,\widehat s)}}\hspace{15mm}\\
&\le& L^{2\over 2\alpha+1}\tau^{2\alpha\over 2\alpha+1}+\sum_{j=1}^{k}
\pa{L(R_{j}k^{-1/2})^{\alpha\wedge 1}}^{2\over 2(\alpha\wedge1)\beta_{j}+1}
\pa{\tau\L}^{2(\alpha\wedge 1)\beta_{j}\over 2(\alpha\wedge 1)\beta_{j}+1}+\tau\L,
\end{eqnarray*}
where $\L=\log\pa{\tau^{-1}}\vee \log\pa{L^{2}}\vee 1$ and $C'$ is a constant that depends 
on $\alpha,\gbeta,\gp$ and $k$ only. 
\end{theorem}
If one is mainly interested in the rate of convergence as $\tau$ tends to 0, the bound we get is of 
order $\max\{\tau^{2\alpha/(2\alpha+1)},[\tau\log(\tau^{-1})]^{2(\alpha\wedge
1)\beta/(2(\alpha\wedge 1)\beta+1)}\}$ where $\beta=\min\{\beta_{1},\ldots,\beta_{k}\}$. In
particular, if $\alpha\ge 1$, this rate is the same as that we would obtain for estimating a function 
on $[-1,1]$ with the smallest regularity among $\alpha,\beta_{1},\ldots,\beta_{k}$.\vspace{2mm}\\
%
{\em Proof:}
Let  us consider some $s=g\circ u\in\S_{k}^{\rm Add}(\alpha,\gbeta,\gp,L,R)$ and $r= 
1+\lceil\alpha\vee \beta_{1}\vee\ldots\vee \beta_{k}\rceil$. For all $D,D_{1},\ldots,D_{k}\ge 1$, 
there exist $F\in\HH_{1,r}(D)$ and $T_{j}\in\BB_{1,r}(D_{j})$ for all $j=1,\ldots,k$ such that 
\[
{\cal D}(F)\le C'_1(r)D;\qquad d_{\infty}(g,F)\le C_2'(r)LD^{-\alpha};
\]
and, for $1\le j\le k$,
\[
{\cal D}(T_j)\le C'_3\!\left(k,r,\gp\right)D_{j};\qquad d(u_j,T_j)\le
C'_4\!\left(k,r,\gp\right)R_{j}k^{-1}D_{j}^{-\beta_j}.
\]
If $T=T_{1}+\ldots+T_{k}$, then $\D(T)\le \sum_{j=1}^{k}\D(T_{j})$,  $\Delta_{\lambda^{(k)}_{r}}
(T)\le \sum_{j=1}^{k}\Delta_{\lambda_r}(T_{j})\le (c(\BB_{1,r})+1)\sum_{j=1}^{k}(D_{j}+1)$.
Moreover, $d(u,T)\le \sum_{j=1}^{k}d(u_{j},T_{k})\le
C_{4}'k^{-1}\sum_{j=1}^{k}R_{j}D_{j}^{-\beta_j}$, hence, $d^{2}(u,T)\le
(C_{4}')^{2}k^{-1}\sum_{j=1}^{k}R_{j}^{2}D_{j}^{-2\beta_j}$ and finally,                
\[
d^{2(\alpha\wedge 1)}(u,T)\le (C_{4}')^{2(\alpha\wedge 1)}
\sum_{j=1}^{k}(R_{j}k^{-1/2})^{2(\alpha\wedge 1)}D_{j}^{-2(\alpha\wedge 1)\beta_j}.
\]
For all $T$, $\L_{1,T}\le C'(\alpha)\L$ and since $\w_{g}(z)\le Lz^{\alpha}$ for all $z\in[0,2]$, 
we may apply Corollary~\ref{C-Main0} with $l=1$ and get that the risk of the resulting estimator
$\widehat s_{r}$ satisfies 
\[
C'\RR(\widehat s_{r},g,u)= \sum_{j=1}^k\inf_{D\ge 1}\cro{L^{2}(R_{j}k^{-1/2})^{2(\alpha
\wedge 1)}D^{-2(\alpha\wedge 1)\beta_{j}}+D\tau\L}+ \inf_{D\ge
1}\cro{L^{2}D^{-2\alpha}+D\tau}.
\]
We conclude by arguing as in the proof of Theorem~\ref{T-Smooth}.\cqfd

\subsection{Multiple index models and artificial neural networks\labs{C3}}
In this section, we assume that $E=[-1,1]^k$, $q=2$ and $d$ is the distance in $\IL_{2}(E,\mu)$ where $\mu$ is the Lebesgue probability on $E$. We denote by $|\cdot|_1$ and
$|\cdot|_\infty$ respectively the $\ell_1$- and $\ell_\infty$-norms in $\R^k$ and
$\C_k$ the unit ball for the $\ell_1$-norm. As we noticed earlier, when $s$ is an arbitrary
function on $E$ and $k$ is large, there is no hope to get a nice estimator for $s$ without some
additional assumptions. A very simple one is that $s(x)$ can be written as $g(\scal{\theta}{x})$
for some $\theta\in\C_k$, which corresponds to the so-called {\em single index model}. More
generally, we may pretend that $s$ can be well approximated by some function $\overline s$ of 
the
form 
\[
\overline s(x)=g\pa{\st\<\theta_{1},x\>,\ldots,\<\theta_l,x\>}
\]
where
$\theta_1,\ldots,\theta_l$ are $l$ elements of $\C_k$ and $g$ maps $[-1,1]^l$  to $\R$,
$l$ being possibly unknown and larger than $k$. When $\bar{s}=g\circ u$ is of this form, the
coordinate functions $u_j(\cdot)=\<\theta_j,\cdot\>$, for $1\le j\le l$, belong to the set
$T_0\subset\T$ of functions on $E$ of the form $x\mapsto\<\theta,x\>$ with
$\theta\in\C_k$, which is a subset of a $k$-dimensional linear subspace of $\IL_2(E,\mu)$, hence
$\D(T_0)\le k$. A slight generalization of this situation leads to the following result.
%
\begin{theorem}\lab{Th-mim}
For $j\ge 1$, let $T_{j}$ be a subset of $\T$ with finite dimension $k_{j}$ and for $I\subset \N^{*}$ and $l\in  I$, let $\FF_{l}$ be a collection of models satisfying Assumptions~\ref{A-models}-(i and iii) for some subprobability $\gamma_{l}$. 
There exists an estimator $\widehat s$ which  satisfies 
\begin{eqnarray}
\lefteqn{C\E_{s}\cro{d^{2}(s,\widehat s)}}\nonumber\hspace{15mm}\\
&\le& \inf_{l\in I}\inf_{g\in \F_{l,c},\,u\in \gT_{l}}\cro{d^2(s,g\circ u)+
A(g,\FF_{l},\gamma_{l})+\tau \sum_{j=1}^{l}k_{j}i(g,j,T_{j})},\qquad
\label{Eq-MainRZ}
\end{eqnarray}
where $\gT_{l}=T_{1}\times\ldots\times T_{l}$, $i(g,j,T_{j})$ is defined by~\eref{Eq-igj} and 
\[
A(g,\FF_{l},\gamma_{l})=\inf_{F\in\FF_{l}}
\ac{d_{\infty}^{2}(g,F)+\tau\cro{\D(F)+\Delta_{\gamma_{l}}(F)}}.
\]
In particular, for all $l\in I$ and $(\galpha,\gL)$-H\"olderian functions $g$ with $\galpha\in(0,1]^{l}$ and $\gL\in(\R_{+}^{\star})^{l}$ 
\begin{eqnarray}
\lefteqn{C\Bbb{E}_s\!\left[d^2\left(s,\widehat s)\right)\st\right]}\nonumber\\
&\le& \inf_{u\in \gT_{l}}\cro{d^2(s,g\circ u)+A(g,\FF_{l},\gamma_{l})+
\tau\sum_{j=1}^lk_{j}\left[{1\over \alpha_{j}}\log\!\left(lL_j^2(k_j\tau)^{-1}\right)
\bigvee 1\right]}.\quad
\labe{Eq-MainRY}
\end{eqnarray}
\end{theorem}
Let us comment on this result, fixing some value $l\in I$. The term $d(s,g\circ u)$ corresponds 
to the approximation of $s$ by functions of the form $g(u_{1}(.),\ldots,u_{l}(.))$ with $g$ in $
\F_{l,c}$ and $u_{1},\ldots,u_{l}$ in $T_{1},\ldots,T_{l}$ respectively. As to the quantity
$A(g,\FF_{l},\gamma_{l})$, it corresponds to the estimation bound for estimating the function
$g$ alone if $s$ were really of the previous form. Finally, the quantity
$\tau\sum_{j=1}^{l}k_{j}i(g,j,T_{j})$ corresponds to the sum of the statistical errors for
estimating the $u_{j}$. If for all $j$, the dimensions of the $T_{j}$ remain bounded by some
integer $\overline k$ independent of $\tau$, which amounts to making a parametric assumption
on the $u_{j}$, and if $g$ is smooth enough the quantity $\tau\sum_{j=1}^{l}k_{j}i(g,j,T_{j})$ is
then of order $\tau\log\tau^{-1}$ for small values of $\tau$ as seen in~\eref{Eq-MainRY}.
\vspace{2mm}\\
%
\noindent{\em Proof of Theorem~\ref{Th-mim}:}
For all $j$, we choose $\lambda_{j}$ to be the Dirac mass at $T_{j}$ so that $\Delta_{\lambda_j}(T_j)=0=d(u_j,T_j)$. 
The result follows by applying Theorem~\ref{T-general} (for a fixed value of $l\in I$) and then
Theorem~\ref{T-melange} with $\nu$ defined by $\nu(l)=e^{-l}$ for all $l\in I$.\cqfd
%

\subsubsection{The multiple index model\labs{C3b}}
As already mentioned,  the multiple index model amounts to assuming that $s$ is of the form 
\[
s(x)=g\pa{\st\<\theta_{1},x\>,\ldots,\<\theta_l,x\>}\quad\mbox{whatever }x\in E,
\]
for some known $l\ge 1$ and $k_j=k$ for all $j$. For $L>0$ and $\galpha\in
(\R_{+}^{\star})^{l}$, let us denote by $\S_{l}^{\galpha}(L)$ the set of functions $s$ of this form
with $g\in \H^{\galpha}([-1,1]^{l})$ satisfying $\norm{g}_{\galpha,\infty}\le L$. Applying
Theorem~\ref{Th-mim} to this special case, we obtain the following result. 
\begin{corollary}
Let $I\subset \N^{\star}$. There exists an estimator $\widehat s$ such that for all $l\in I$, 
$\galpha\in(\R_{+}^{\star})^{l}$ and $L>0$, 
\begin{eqnarray*}
\sup_{s\in\S_{l}^{\galpha}(L)}C'\Bbb{E}_s\!\left[d^2\left(s,\widehat s\right)\st\right]&
\le& L^{{2\over 2\overline{\alpha}+1}}\tau^{2\overline{\alpha}\over
2\overline{\alpha}+1}+k\tau\L,
\end{eqnarray*}
where $\L=\log(\tau^{-1})\vee \log (L^{2}k^{-1})\vee 1$ and $C'$ is a constant depending on 
$l$ and $\galpha$ only. 
\end{corollary}
\noindent{\em Proof:}
Fix $s=g\circ u\in\S_{l}^{\galpha}(L)$ and apply Theorem~\ref{Th-mim} with $T_{j}=T_{0}$
for all $j\ge 1$, $I=\{l\}$, $\FF_{l}=\HH_{l,r}$ and $\gamma_{l}$ defined by~\eref{Eq-defgamma} 
with $k=l$ and $r=\lceil \alpha_{1}\vee\ldots\vee \alpha_{l}\rceil$. Arguing as in the proof of
Theorem~\ref{T-Smooth}, we obtain an estimator $\widehat s_{(l,r)}$ the risk of which satisfies 
\[
\RR(\widehat s_{(l,r)},g,u)=C'\cro{\inf_{D\ge 1}\pa{L^{2}D^{-2\overline \alpha/l}+D
\tau}+ \tau k\L}\le C''\cro{L^{{2\over 2\overline{\alpha}+1}}\tau^{2\overline{\alpha}\over
2\overline{\alpha}+1}+k\tau\L},
\]
for constants $C'$ and $C''$ depending on $l$ and $\galpha$ only. Finally, we conclude as in 
the proof of Theorem~\ref{T-Smooth}.\cqfd

\subsubsection{Case of an additive function $g$\labs{C3a}}
In the multiple index model, when the value of $l$ is allowed to become large (typically not 
smaller than $k$) it is often assumed that $g$ is additive, i.e.\ of the form 
\begin{equation}
g(y_{1},\ldots,y_{l})=g_{1}(y_{1})+\ldots+ g_{l}(y_{l})\quad\mbox{for all }y\in[-1,1]^{l},
\labe{Eq-add}
\end{equation}
where the $g_j$ are smooth functions from $[-1,1]$ to $\R$. Hereafter, we shall denote by $\F_{l,c}^{{\rm Add}}$ the set of such additive functions $g$. The functions  $\overline s=g\circ u$ with $g\in \F_{l,c}^{{\rm Add}}$ and $u\in\ T_{0}^{l}$ hence  take the form
\begin{equation}
\overline s(x)=\sum_{j=1}^lg_j\pa{\<\theta_j,x\>}\quad\mbox{for all }x\in E.
\label{add-s}
\end{equation}
For each $j=1,\ldots,l$, let $\FF_j$ be a countable family of finite dimensional linear subspaces
of $\F_{1,\infty}$ designed to approximate $g_j$ and
$\gamma_j$ some subprobability measure on $\FF_j$. Given $(F_1,\ldots,F_l)\in
\prod_{j=1}^l\FF_j$, we define the subspace $F$ of $\F_{l,\infty}$ as
\begin{equation}
F=\ac{f(y_{1},\ldots,y_{l})=f_{1}(y_{1})+\ldots+ f_{l}(y_{l})\,\left|\,f_j\in F_j 
\mbox{ for } 1\le j\le l\st\right.}
\labe{Eq-FF}
\end{equation}
and denote by $\FF$ the set of all such $F$ when $(F_1,\ldots,F_l)$ varies among
$\prod_{j=1}^l\FF_j$. Then, we define a subprobability measure $\gamma$ on $\FF$ by
setting
\[
\gamma(F)=\prod_{j=1}^l\gamma_j(F_j)\quad\mbox{or}\quad\Delta_{\gamma}(F)=
\sum_{j=1}^l\Delta_{\gamma_j}(F_j),
\]
when $F$ is given by (\ref{Eq-FF}). For such an $F$, $d_{\infty}(g,F)\le\sum_{j=1}^l
d_{\infty}(g_j,F_j)$, hence $d^2_{\infty}(g,F)\le l\sum_{j=1}^ld^2_{\infty}(g_j,F_j)$ and $\D(F)\le \sum_{j=1}^l\D(F_j)$. We deduce from Theorem~\ref{Th-mim} the following result. 
\begin{corollary}\lab{C-add}
Let $I\subset \N^{\star}$ and for $j\ge 1$, let $\FF_{j}$ be a collection of finite dimensional 
linear subspaces of $\F_{1,\infty}$  satisfying Assumption~\ref{A-models}-i) and-iii)  for some
subprobability $\gamma_{j}$. There exists an estimator $\widehat s$ such that 
\[
C\Bbb{E}\cro{d^{2}(s,\widehat s)}\le \inf_{l\in I}\inf_{g\in\F_{l,c}^{{\rm Add}},u\in 
T_{0}^{l}}\cro{d^{2}(s,g\circ u)+\sum_{j=1}^{l}\pa{\st
R_{j}(g,\FF_{j},\Delta_{\gamma_{j}})+\tau ki(g,j,T_{0})}},
\]
where
\[
R_j(g,\FF_{j},\Delta_{\gamma_{j}})=\inf_{F_j\in\FF_j}\ac{d_{\infty}^2(g_j,F_j)+\tau
\cro{\D(F_j)+\Delta_{\gamma_j}(F_j)}\st}\quad\mbox{for }1\le j\le l.
\]
Moreover, if $s$ of the form~\eref{add-s} for some $l\in I$ and functions $g_{j}\in
\H^{\alpha_{j}}([-1,1])$ satisfying $\norm{g_{j}}_{\alpha_{j},\infty}\le L_{j}$ for 
$\alpha_{j},L_{j}>0$ and all $j=1,\ldots,l$, one can choose the $\FF_{j}$ and $\gamma_{j}$ in 
such a way that
\begin{equation}\label{cash}
\E_{s}\cro{d^{2}(s,\widehat s}\le C'\cro{\sum_{j=1}^{l}L_{j}^{2\over 2\alpha_{j}+1}
\tau^{2\alpha_{j}\over 2\alpha_{j}+1}\ +\ k\tau\L},
\end{equation}
where $\L=\log\pa{\tau^{-1}}\vee 1\vee\cro{\bigvee_{j=1}^{l}\log\pa{L_{j}^{2}k^{-1}}}$ and 
$C'$ is a constant depending on $l$ and $\alpha_{1},\ldots,\alpha_{l}$ only. 
\end{corollary}

For $j\ge1$, $R_{j}=R_j(g,\FF_{j},\Delta_{\gamma_{j}})$ corresponds to the risk bound for the 
estimation of the function $g_j$ alone when we use the family of models $\FF_j$, i.e.\ what we
would get if we knew $\theta_j$ and that $g_i=0$ for all $i\ne j$. In short, $\sum_{j=1}^lR_j$
corresponds to the estimation rate of the additive function $g$. If each $g_j$ belongs to some
smoothness class, this rate is similar to that of a real-valued function defined on the line with
smoothness given by the worst component of $g$, as seen in~\eref{cash}.\vspace{2mm}

%
\noindent{\em Proof of Corollary~\ref{C-add}:}
The first part is a straightforward consequence of Theorem~\ref{Th-mim}. For the second part, 
fix $s=g\circ u$ and $r=\lfloor \alpha_{1}\vee\ldots\vee \alpha_{l}\rfloor$. Since the $g_{j}$
are $(\alpha_{j}\wedge 1,L_{j})$-H\"olderian, $i(g,j,T_{0})\le C'\L$ for some $C'$ depending on
the $\alpha_{j}$ only. By using Proposition~\ref{P-Approx1}, Lemma~\ref{L-opt} and the
collection $\FF_{j,r}=\HH_{1,r}$ with $\gamma_{j,r}$ defined by~\eref{Eq-defgamma}, for all
$j=1,\ldots,l$, $R_{j}\le C'\inf_{D\ge 1}\{L_{j}^{2}D^{-2\alpha_{j}}+D\tau\}\le
C''(L_{j}^{2/(2\alpha_{j}+1)}\tau^{2\alpha_{j}/(2\alpha_{j}+1)}+\tau)$ for some constants
$C',C''$ depending on the $\alpha_{j}$ only. Putting these bounds together, we end up with an
estimator $\widehat s_{r}$ the risk of which is bounded from above by the right-hand side
of~\eref{cash}. We get the result for all values of $r$ by using Theorem~\ref{T-melange} and
arguing as in the proof of Theorem~\ref{T-Smooth}. \cqfd
%
%
\subsubsection{Artificial neural networks\labs{C3c}}
In this section, we consider approximations of $s$ on $E=[-1,1]^k$ by functions of the form
\begin{equation}
\bar{s}(x)=\sum_{j=1}^lR_j\psi\pa{\<a_j,x\>+b_j}\quad\mbox{with }|b_j|+|a_j|_1\le2^q,
\labe{Eq-ANN}
\end{equation}
for given values of $(l,q)\in I=(\N^{\star})^2$. Here, $R=(R_1,\ldots,R_l)\in\R^l$, $a_j\in
\R^k$, $b_j\in\R$ for $j=1,\ldots,l$ and $\psi$ is a given uniformly continuous function on
$\Bbb{R}$ with modulus of continuity  $\w_{\psi}$. We denote by $\S_{l,q}$ the set of all
functions $\bar{s}$ of the form (\ref{Eq-ANN}).

Let us  now set $\psi_q(y)=\psi\left(2^qy\right)$ for $y\in\Bbb{R}$ and, for $x\in E$,
$u_j(x)=2^{-q}\pa{\<a_j,x\>+b_j}$, so that $u_j\in\T$ belongs to the $(k+1)$-dimensional 
spaces of functions of the form $x\mapsto\<a,x\>+b$. We can then rewrite
$\bar{s}$ in the form $g\circ u$ with $g(y_1,\ldots,y_l)=\sum_{j=1}^l R_j\psi_q(y_j)$.
Since $g$ belongs to the $l$-dimensional linear space $F$ spanned by the functions
$\psi_q(y_j)$, we may set $\FF=\{F\}$, $\Delta_\gamma(F)=0$ and apply
Theorem~\ref{Th-mim}. With $\w_{g,j}(y)=|R_j|\w_\psi\left(2^qy\right)$,
(\ref{Eq-MainRZ}) becomes, 
\[
C\Bbb{E}_s\!\left[d^2\left(s,\widehat{s}_{l,q}\right)\st\right]\le d^2(s,\bar{s})+
\tau(k+1)\sum_{j=1}^l
\inf\left\{i\in\N^\star\,\left|\,lR_j^2\w_\psi^2\left(2^qe^{-i}\right)\le(k+1)\tau
i\right.\right\}.
\]
If $\w_\psi(y)\le Ly^\alpha$ for some $L>0$, $0<\alpha\le1$ and all $y\in\R_+$, then, 
according to (\ref{Eq-LJT}),
\begin{eqnarray}
C\Bbb{E}_s\!\left[d^2\left(s,\widehat{s}_{l,q}\right)\st\right]
\nonumber&\le&d^2(s,\bar{s})+k\tau\left(\sum_{j=1}^l
\left[\alpha^{-1}\log\!\left(lR_j^2L^22^{2q\alpha}[k\tau]^{-1}\right)\right]
\bigvee1\right)\nonumber\\&\le&d^2(s,\bar{s})+lk\tau\left[q\log4+\alpha^{-1}
\log_+\!\left(l\left|R\right|_\infty^2L^2[k\tau]^{-1}\right)\right].\qquad\quad
\labe{Eq-ANN1}
\end{eqnarray}
These bounds being valid for all $(l,q)\in I$ and $\bar{s}\in\S_{l,q}$, we may apply
Theorem~\ref{T-melange} to the family of all estimators $\widehat{s}_{l,q},\,(l,q)\in I$ with
$\nu$ given by $\nu(l,q)=e^{-l-q}$. We then get the following result.
%
\begin{theorem}\lab{T-ANN}
Assume that $\psi$ is a continuous function with modulus of continuity $\w_\psi(y)$
bounded by  $Ly^\alpha$ for some $L>0$, $0<\alpha\le1$ and all $y\in\R_+$. Then one 
can build an estimator $\widehat{s}=\widehat{s}({\bf X})$ such that 
\begin{eqnarray}
\lefteqn{C\Bbb{E}_s\!\left[d^2\left(s,\widehat{s})\right)\st\right]}\hspace{1mm}
\nonumber\\&\le&\inf_{(l,q)\in I}\,\inf_{\bar{s}\in\S_{l,q}}\left\{d^2(s,\bar{s})+lk\tau q
\left[1+(q\alpha)^{-1}\log_+\!\left(l\left|R\right|_\infty^2L^2[k\tau]^{-1}\right)\right]
\right\}.\qquad\quad
\labe{Eq-ANN2}
\end{eqnarray}
\end{theorem}
%

\paragraph{Approximation by functions of the form (\ref{Eq-ANN}).}
Various authors have provided conditions on the function $s$ so that it can be approximated 
within $\eta$ by functions $\bar{s}$ of the form (\ref{Eq-ANN}) for a given function $\psi$. An
extensive list of authors and results is provided in Section~4.2.2 of Barron, Birg\'e and Massart~\citeyearpar{MR1679028} and some proofs are provided in Section~8.2 of that paper. The starting point of such
approximations is the assumed existence of a Fourier representation of $s$ of the form
\[
s(x)=K_s\int_{\Bbb{R}^k}\cos\left(\st\scal{a}{x}+\delta(a)\right)dF_s(a),\quad
K_s\in\Bbb{R},\quad|\delta(a)|\le\pi,
\]
for some {\em probability} measure $F_s$ on $\Bbb{R}^k$. To each given function $\psi$ that can
be used for the approximation of $s$ is associated a positive number $\beta=\beta(\psi)>0$ and
one has to assume that 
\begin{equation}
c_{s,\beta}=\int|a|^\beta_1dF_s(a)<+\infty,
\labe{Eq-ANN7}
\end{equation}
in order to control the approximation of $s$ by functions of the form (\ref{Eq-ANN}). 
A careful inspection of the proof of Proposition~6 in Barron, Birg\'e and Massart ~\citeyearpar{MR1679028} shows that,
when (\ref{Eq-ANN7}) holds, one can derive the following approximation result for $s$. There
exist constants $q_\psi\ge1$, $\gamma_\psi>0$ and $C_\psi>0$ depending on $\psi$ only, a
number $R_{s,\beta}\ge1$ depending on $c_{s,\beta}$ only and some $\bar{s}\in\S_{l,q}$ with
$|R|_1\le R_{s,\beta}$ such that
\begin{equation}
d\left(s,\bar{s}\right)\le K_sC_\psi\left[2^{-q\gamma_\psi}+
R_{s,\beta}l^{-1/2}\right]\quad\mbox{for } q\ge q_\psi.
\labe{Eq-ANN3}
\end{equation}
Putting this bound into (\ref{Eq-ANN2}) and omitting the various indices for simplicity, we
get a risk bound of the form
\[
\RR(l,q)=CK^2\left[2^{-2q\gamma}+R^2l^{-1}+K^{-2}lk\tau q\left[1+
(q\alpha)^{-1}\log_+\left(lR^2L^2[k\tau]^{-1}\right)\right]\right],
\]
to be optimized with respect to $l\ge1$ and $q\ge q_\psi$. We shall actually perform the optimization with respect to the first three terms, omitting the logarithmic one.

Let us first note that, if $RK\le\sqrt{q_\psi k\tau}$, one should set $q=q_\psi$ and $l=1$, which leads to
\[
\RR(1,q_\psi)\le Ck\tau q_\psi\left[1+(q_\psi\alpha)^{-1}\log_+\left(R^2L^2[k\tau]^{-1}\right)\right].
\]
Otherwise $\sqrt{q_\psi k\tau}<RK$ and we set
\[
q=q^{*}=\inf\ac{q\ge q_{\psi}\ \Big|\ 2^{-2q\gamma}\le(R/K)\sqrt{qk\tau}}\qquad\mbox{and}\qquad
l=l^*=\left\lceil\frac{RK}{\sqrt{q^*k\tau}}\right\rceil.
\]
If $l^*>1$, then $RK(q^*k\tau)^{-1/2}\le l^*<2RK(q^*k\tau)^{-1/2}$ hence
\begin{equation}
\RR(l^*,q^*)\le CRK\sqrt{q^*k\tau}\left[1+\frac{1}{q^*\alpha}\log_+\left(\frac{2R^3L^2K}{(k\tau)^{3/2}\sqrt{q^*}}\right)\right].
\labe{Eq-ANN4}
\end{equation}
If $l^*=1$, then $R^2\le K^{-2}q^*k\tau$ and $\sqrt{q_\psi k\tau}<RK\le\sqrt{q^*k\tau}$, hence $q^*>q_\psi$ and $q^*-1\ge q^*/2$. Then, from the definition of $q^*$,
\[
RK^{-1}\sqrt{(q^*/2)k\tau}\le RK^{-1}\sqrt{(q^*-1)k\tau}<2^{-2(q^*-1)\gamma}\le2^{-2\gamma},
\]
hence $\sqrt{q^*k\tau}<(K/R)2^{-2\gamma+(1/2)}<\sqrt{2}K$ and (\ref{Eq-ANN4}) still holds. To conclude, we observe that either $-2\gamma q_\psi\log2\le\log\left(RK^{-1}\sqrt{q_\psi k\tau}\right)$ and $q^*=q_\psi$ or the solution $z_0$ of the equation
\[
2z\gamma\log2=\log\left(K/\left[R\sqrt{k\tau}\right]\right)-(1/2)\log z
\]
satisfies $q_\psi<z_0\le q^*$. Since $\log z_0\le z_0/e$, it follows that
\[
q^*\ge\log\left(K/\left[R\sqrt{k\tau}\right]\right)/\left(2\gamma\log2+e^{-1}\right)
\]
and, by monotonicity, that
\[
\frac{1}{q^*}\log_+\left(\frac{2R^3L^2K}{(k\tau)^{3/2}\sqrt{q^*}}\right)\le\L=
\left(2\gamma\log2+e^{-1}\right)\log_+\left(\frac{2R^3L^2K}{(k\tau)^{3/2}\sqrt{q_\psi}}\right)
\left[\log\left(\frac{K}{R\sqrt{k\tau}}\right)\right]^{-1}
\]
where $\L$ is a bounded function of $k\tau$. One can also check that
\[
q^*\le\overline q=\left\lceil\frac{\log\left(K/\left[R\sqrt{q_\psi k\tau}\right]\right)}{2\gamma\log2}\right\rceil
\]
and (\ref{Eq-ANN4}) finally leads, when $q^*>q_\psi$, to
\begin{equation}
\RR(l^*,q^*)\le CRK\left(k\tau\left\lceil\frac{\log\left(K/\left[R\sqrt{q_\psi k\tau}\right]\right)}{2\gamma\log2}\right\rceil\right)^{1/2}\left[1+\alpha^{-1}\L\right].
\labe{Eq-ANN5}
\end{equation}
In the asymptotic situation where $\tau$ converges to zero, (\ref{Eq-ANN5}) prevails and we get a risk bound of order $\cro{-k\tau\log(k\tau)}^{1/2}$.

\subsection{Estimation of a regression function and PCA\labs{C2}}
We consider here the regression framework 
\[
Y_{i}=s(X_{i})+\eps_i,\quad i=1,\ldots,n,
\]
where the $X_{i}$ are random variables with values in some known compact subset $K$ of
$\R^{k}$ (with $k>1$ to avoid trivialities) the $\eps_{i}$ are i.i.d.\ centered random variables
of common variance 1 for simplicity  and $s$ is an unknown function from $\R^{k}$ to 
$\R$. By a proper origin and scale change on the $X_i$, mapping $K$ into the unit ball $\B_k$
of $\R^{k}$, one may assume that the $X_i$ belong to $\B_k$, hence that $E=\B_k$, which we
shall do from now on. We also assume that the $X_i$ are either i.i.d.\ with common distribution
$\mu$ on $E$ (random design) or deterministic ($X_i=x_i$, fixed design), in which case
$\mu=n^{-1}\sum_{i=1}^n\delta_{x_i}$, where $\delta_x$ denotes the Dirac measure at $x$. In 
both cases, we choose for $d$ the distance in $\mathbb{L}_2(E,\mu)$. As already mentioned in 
Section~\ref{F2}, Theorem~\ref{T-models} with $\tau=n^{-1}$ applies to this framework, at least in the two cases when the design is fixed and the errors Gaussian (or subgaussian) or when the design is random and the $Y_i$ are bounded, say with values in $[-1,1]$. 

\subsubsection{Introducing PCA\labs{C2z}}
Our aim is to estimate $s$ from the observation of the pairs $(X_{i},Y_{i})$ for $i=1,\ldots,n$,
assuming that $s$ belongs to some smoothness class. More precisely, given $A\subset\R^k$ and  some concave modulus of continuity $\w$ on $\R_{+}$, we define $\H_{\w}(A)$ to be the class of functions $h$ on $A$ such that 
\[
\ab{h(x)-h(y)}\le \w\pa{\ab{x-y}}\quad\mbox{for all }x,y\in A.
\]  
Here we assume that $s$ is defined on $\B_k$ and belongs to  $\H_{\w}(\B_k)$, in which case
it can be extended to an element of $\H_{\w}\left(\R^k\right)$, which we shall use when
needed. Typically, if $\w(z)=Lz^{\alpha}$ with $\alpha\in(0,1]$ and the $X_{i}$ are i.i.d.\ with 
uniform distribution $\mu$ on $E$, the minimax risk bound over
$\H_{\w}\left(\B_k\right)$ with respect to the $\mathbb{L}_2(E,\mu)$-loss is
$C'L^{2k/(k+2\alpha)}n^{-2\alpha/(k+2\alpha)}$ (where $C'$ depends on $k$ and the
distribution of the $\varepsilon_i$). It can be quite slow if $k$ is large (see Stone~\citeyearpar{MR673642}),
although no improvement is possible from the minimax point of view  if the distribution of the
$X_i$ is uniform on $\B_k$. Nevertheless, if the data $X_i$ were known to belong to an affine
subspace $V$ of $\R^k$ the dimension $l$ of which is small as compared to
$k$, so that $\mu(V)=1$, estimating the function $s$ with  $\mathbb{L}_2(E,\mu)$-loss  would
amount to estimating $s\circ\Pi_{V}$ (where $\Pi_{V}$ denotes the orthogonal projector onto
$V$) and one would get the much better rate $n^{-2\alpha/(l+2\alpha)}$ with respect to $n$ for
the quadratic risk. Such a situation is seldom encountered in practice but we may assume that it
is approximately satisfied for some well-chosen $V$. It therefore becomes natural to look for an
affine space $V$ with dimension $l<k$ such that $s$ and $s\circ\Pi_{V}$ are close with respect
to the $\mathbb{L}_2(E,\mu)$-distance. For $s\in\H_{\w}\left(\R^k\right)$, it follows from
Lemma~\ref{L-approx} below that, 
\begin{eqnarray*}
\int_{E}\ab{s(x)-s\circ\Pi_{V}(x)}^{2}d\mu(x)&\le& \int_{E}\w^{2}\pa{\ab{x-\Pi_{V}x}}
d\mu(x)\\ &\le&2\w^{2}\cro{\pa{\int_{E}\ab{x-\Pi_{V}x}^{2}d\mu(x)}^{1/2}},
\end{eqnarray*}
and minimizing the right-hand side amounts to finding an affine space $V$ with dimension 
$l$ for which $\int_{E}\ab{x-\Pi_{V}x}^{2}d\mu(x)$ is minimum. This way of reducing the
dimension is usually known as PCA (for Principal  Components Analysis). When the $X_{i}$
are deterministic and $\mu=n^{-1}\sum_{i=1}^{n}\delta_{X_{i}}$, the solution to this
minimization problem is given by the affine space $V_l=a+W_l$ where the origin
$a=\overline{X}_{n}=n^{-1}\sum_{i=1}^{n}X_{i}\in\B_k$ and $W_l$ is the linear space
generated by the eigenvectors associated to the $l$ largest eigenvalues (counted with their
multiplicity) of  $XX^{*}$ (where $X$ is the  $k\times n$  matrix with columns
$X_i-\overline{X}_{n}$ and $X^*$ is the transpose of $X$). In the general case, it suffices to set
$a=\int_{E} xd\mu$ (so that $a\in E$) and replace $XX^{*}$ by the matrix 
\[
\Gamma= \int_{E}(x-a)(x-a)^{*}\,d\mu(x).
\]
If $\lambda_{1}\ge \lambda_{2}\ge \ldots\ge \lambda_{k}\ge 0$ are the eigenvalues of 
$\Gamma$ in nonincreasing order, then 
\begin{equation}
\inf_{\{V\,|\,\dim(V)=l\}}\int_{E}\ab{x-\Pi_{V}x}^{2}d\mu(x)=\sum_{j=l+1}^k\lambda_j
\labe{Eq-rigolo}
\end{equation}
(with the convention $\sum_{\varnothing}=0$) and therefore
\begin{equation}
\inf_{\{V\,|\,\dim(V)=l\}}\norm{s-s\circ\Pi_{V}}_{2}^{2}\le\norm{s-s\circ\Pi_{V_l}}_{2}^{2}
\le2\w^{2}\pa{\sqrt{\sum_{j=l+1}^k\lambda_j}}.
\labe{Eq-rigo}
\end{equation}

\subsubsection{PCA and composite functions\labs{C2x}}
In order to put the problem at hand into our framework, we have to express $s\circ\Pi_{V_l}$
in the form $g\circ u$. To do so we consider an orthonormal basis $\overline{u}_{1},\ldots,
\overline{u}_{k}$ of eigenvectors of  $XX^{*}$ or $\Gamma$ (according to the situation)
corresponding to the  ordered eigenvalues $\lambda_1\ge\lambda_2\ge\ldots\ge
\lambda_k\ge0$. For a given value of $l<k$ we denote by  $a^\perp$ the component of $a$
which is orthogonal to the linear span $W_l$ of $\overline{u}_{1},\ldots,\overline{u}_l$ and
for $x\in\B_k$, we define $u_j(x)=\<x,\overline{u}_j\>$ for $j=1,\ldots,l$. This results in an
element $u=(u_1,\ldots,u_l)$ of $\T^{l}$ and $a^\perp+\sum_{j=1}^lu_j(x)\overline{u}_j
=\Pi_{V_l}(x)$ is the projection of $x$ onto the affine space $V_l=a^\perp+W_l$. Setting
\[
g(z)=s\left(a^\perp+\sum_{j=1}^lz_j\overline{u}_j\right)\quad\mbox{for }z\in[-1,1]^{l},
\]
leads to a function $g\circ u$ with $u\in\T^{l}$ and $g\in\F_{l,c}$ which coincides with
$s\circ\Pi_{V_l}$ on $\B_k$ as required. Consequently, the right-hand side of (\ref{Eq-rigo})
provides a bound on the distance between $s$ and  $g\circ u$. Moreover, since $s\in
\H_{\w}\left(\R^k\right)$,
\begin{equation}
\ab{g(z)-g(z')}\le\w\pa{\ab{\sum_{j=1}^lz_j\overline{u}_j-\sum_{j=1}^lz'_j\overline{u}_j}}
=\w\pa{\ab{\sum_{j=1}^l(z_j-z'_j)\overline{u}_j}}=\w(|z-z'|),
\labe{Eq-RB7}
\end{equation}
so that we may set $\w_{g,j}=\w$ for all $j\in\ac{1,\ldots,l}$.

In the following sections we shall use this  preliminary result in order to establish risk bounds for
estimators $\widehat{s}_l$ of $s$, distinguishing between the two situations where $\mu$ is
known and $\mu$ is unknown. 

\subsubsection{Case of a known $\mu$\labs{C2a}}
For $D\in\Bbb{N}^\star$, we consider the partition $\PP_{l,D}$ of  $[-1,1]^{l}$ into $D^l$
cubes with edge length $2/D$ and denote by $F_{l,D}$ the linear space of functions which are
piecewise constant on each element of $\PP_{l,D}$ so that $\D(F_{l,D})=D^l$ for all
$D\in\Bbb{N}^\star$. This leads to the family $\FF=\{F_{l,D},\,D\in\Bbb{N}^\star\}$ and
we set $\gamma(F_{l,D})= e^{-D}$ for all $D\ge 1$. We define $u_j$ as in the previous section
and take for $\TT_j$ the family reduced to the single model $T_j=\ac{u_j}$ for $j=1,\ldots,l$.
Then $\D(T_j)=0$ for all $j$ and we take for $\lambda_j$ the Dirac measure on $\TT_j$. This
leads to a set $\frak{S}$ which satisfies Assumption~\ref{A-models} and we may therefore
apply Theorem~\ref{T-general} which leads to an estimator
$\widehat{s}_l$ with a risk bounded by
\[
C\Bbb{E}_s\cro{\norm{s-\widehat{s}_l}_{2}^{2}}\le d^2(s,g\circ u)+
\inf_{D\ge1}\ac{d_{\infty}^{2}(g,F_{l,D})+{D^l+D\over n}}.
\]
Since $s\circ \Pi_{V_l}$ and $g\circ u$ coincide on $\B_k$, it follows from \eref{Eq-rigo} that
\[
\norm{s-g\circ u}_{2}^{2}=\norm{s-s\circ \Pi_{V_l}}_{2}^{2}\le
2\w^{2}\pa{\sqrt{\sum_{j=l+1}^k\lambda_j}}.
\]
Moreover, for all cubes $I\in\PP_{l,D}$ and $x\in I$, the Euclidean distance between $x$ and 
the center of $I$ is at most $\sqrt{l}D^{-1}$, hence by (\ref{Eq-RB7}),  $d_{\infty}(g,F_{l,D})\le
\w\left(\sqrt{l}D^{-1}\right)$ for all $D\ge 1$. Putting these inequalities together we see that
the risk of
$\widehat{s}_l$ is bounded by
\begin{equation}
C\Bbb{E}_s\cro{\norm{s-\widehat{s}_l}_{2}^{2}}\le \w^2\pa{\sqrt{\sum_{j=l}^k\lambda_j}}
+\inf_{D\ge1}\left\{\w^2\left(\sqrt{l}D^{-1}\right)+{D^l\over n}\right\}.
\labe{Eq-RB8}
\end{equation}

\subsubsection{Case of an unknown $\mu$\labs{C2b}}
When $\mu$ corresponds to an unknown distribution of the $X_{i}$, the matrix $\Gamma$
is unknown, its eigenvectors $\overline{u}_{1},\ldots, \overline{u}_{k}$ and the vector $a$
as well and therefore also the elements $u_1,\ldots,u_l$ of $\T$. In order to cope with
this problem, we have to approximate the unknown $u_j$ which requires to modify the
definition of $T_j$ given in the previous section, keeping all other things unchanged. For each 
$v\in\R^k$ with $|v|\le1$, we denote by $t_v$ the linear map, element of $\T$, given by
$t_v(x)=\<x,v\>$. Denoting by ${\cal B}^{\circ}_k$ the unit sphere in $\R^k$ we then set, for all
$j$, $T_j=T=\{t_v,\,v\in\B_k^{\circ}\}$ which is a subset of a $k$-dimensional linear subspace of
$\mathbb{L}_2(\mu)$. It follows that Assumption~\ref{A-models} remains satisfied but now with
$\D(T_j)=k$. Since $u_{j}\in T_j$ for all $j$, an application of Theorem~\ref{T-general} leads to
\[
C\Bbb{E}_s\!\left[d^2\left(s,\widehat s\right)\st\right]\le\frac{k}{n}\sum_{j=1}^l i(g,j,T)+
d^2(s,g\circ u)+\inf_{D\ge1}\ac{d_{\infty}^{2}(g,F_{l,D})+{D^l+D\over n}},
\]
where $i(g,j,T)$ is given by (\ref{Eq-igj}). Since, by (\ref{Eq-RB7}), $\w_{g,j}=\w$ for all
$j\in\ac{1,\ldots,l}$, 
\[
i(g,j,T)=\underline{i}=\inf\left\{i\in\N^\star\,|\,l\w^2\left(e^{-i}\right)\le\frac{ik}{n}\right\}.
\]
Arguing as in the case of a known $\mu$, we get
\[
C\Bbb{E}_s\cro{\norm{s-\widehat{s}_l}_{2}^{2}}\le \frac{kl\underline{i}}{n}+
\w^2\pa{\sqrt{\sum_{j=l+1}^k\lambda_j}}+\inf_{D\ge1}\left\{
\w^2\left(\sqrt{l}D^{-1}\right)+{D^l\over n}\right\}.
\]
Let $i_D=\left\lceil\log\pa{D/\sqrt{l}}\right\rceil$. If $\underline{i}\le i_D$, then
$kl\underline{i}/n\le klD/n$ since $i_D\le D$. Otherwise, $\underline{i}\ge i_D+1\ge2$ and
\[
l^2\w^2\left(e^{-i_D}\right)\ge
l^2\w^2\left(e^{-\underline{i}+1}\right)>\frac{kl(\underline{i}-1)}{n}\ge\frac{kl\underline{i}}{2n},
\]
which shows that $kl\underline{i}/n<2l^2\w^2\left(\sqrt{l}D^{-1}\right)+klD/n$. Finally
\begin{eqnarray*}
C\Bbb{E}_s\cro{\norm{s-\widehat{s}_l}_{2}^{2}}&\le&\w^2\pa{\sqrt{\sum_{j=l+1}^k\lambda_j}}+
\inf_{D\ge1}\left\{l^2\w^2\left(\sqrt{l}D^{-1}\right)+{D^l+klD\over n}\right\}\\&\le&
\w^2\pa{\sqrt{\sum_{j=l+1}^k\lambda_j}}+k^2\inf_{D\ge1}\left\{
\w^2\left(\sqrt{l}D^{-1}\right)+{2D^l\over n}\right\},
\end{eqnarray*}
which is, up to constants, the same as (\ref{Eq-RB8}).

\subsubsection{Varying $l$\labs{C2c}}
The previous bounds are valid for all values of $l\in I=\{1,\ldots,k\}$ but we do not know
which value of $l$ will lead to the best estimator. We may therefore apply
Theorem~\ref{T-melange} with $\nu(l)=l^{-2}/2$ for $l\in I$ which leads to the following
risk bound for the new estimator $\widehat{s}$ in the case of a known $\mu$:
\[
C\Bbb{E}_s\cro{\norm{s-\widehat{s}}_2^2}\le\inf_{l\in\ac{1,\ldots,k}}\,\inf_{D\ge1}
\cro{\w^2\pa{\sqrt{\sum_{j=l+1}^k\lambda_j}}+\w^2\left(\sqrt{l}D^{-1}\right)+{D^l+\log l\over n}}.
\]
Apart from multiplicative constants depending only on $k$, the same result holds when $\mu$
is unknown. If $\w(z)=Lz^{\alpha}$ for some $L>0$ and $\alpha\in(0,1]$, we get, since
$\sum_{j=l+1}^k\lambda_j\le(k-l)\lambda_{l+1}$ (with the convention $\lambda_{k+1}=0$),
\[
C\Bbb{E}_s\cro{\norm{s-\widehat{s}}_2^2}\le\inf_{l\in\ac{1,\ldots,k}}\,\inf_{D\ge 1}
\ac{L^2[(k-l)\lambda_{l+1}]^{\alpha}+L^2l^\alpha D^{-2\alpha}+{D^l+\log l\over n}}.
\]
Assuming that $n\ge L^{-2}$ to avoid trivialities and choosing
$D=\left\lfloor\left(nL^2l^\alpha\right)^{1/(l+2\alpha)}\right\rfloor$, we finally get 
\[
C\Bbb{E}_s\cro{\norm{s-\widehat s}_2^2}\le
\inf_{l\in\ac{1,\ldots,k}}\ac{L^2[(k-l)\lambda_{l+1}]^{\alpha}+{\log l\over n}
+\frac{L^{2l/(l+2\alpha)}}{n^{2\alpha/(l+2\alpha)}}}.
\]
For $l=k$, we recover (up to constants) the minimax risk bound over $\H_{\w}(\B_k)$, namely
$C'(k)L^{2k/(k+2\alpha)}n^{-2\alpha/(k+2\alpha)}$. Therefore our procedure can only
improve the risk as compared to the minimax approach.

\subsection{Introducing parametric models}
In this section, we approximate $s$ by functions of the form $\overline s=g\circ u$ where $g$ belongs to $\F_{l,c}$ and the components $u_{j}$ of $u$ to parametric models $\gT_{j}=\{u_{j}(\gtheta,.),\ \gtheta\in \Theta_{j}\}\subset \T$ indexed by subsets $\Theta_{j}$ of $\R^{k_{j}}$ with $k_{j}\ge 1$. Besides, we assume that the following holds. 
\begin{assumption}\label{A-embed}
For each $j=1,\ldots,l$, $\Theta_{j}\subset \B_{k_{j}}(0,M_{j})$ for some positive number $M_{j}$ and the mapping $\gtheta\mapsto u_{j}(\gtheta,.)$ from $\Theta_{j}$ to $(\T,d)$ is $(\beta_{j},R_{j})$-H\"olderian for $\beta_{j}\in(0,1]$ and $R_{j}>0$ which means that
\begin{equation}
d(u_{j}(\gtheta,.),u_{j}(\gtheta',.))\le R_{j}\ab{\gtheta-\gtheta'}^{\beta_{j}}\quad\mbox{for all}\quad
\gtheta,\gtheta'\in\Theta_{j}.
\labe{Eq-embed}
\end{equation}
\end{assumption}
Under such an assumption, the following result holds. 

\begin{theorem}\lab{T-para}
Let $l\ge 1$, $\gT_{1},\ldots,\gT_{l}$ be parametric sets satisfying Assumption~\ref{A-embed}, 
$\FF$  be a collection of models satisfying Assumption~\ref{A-models}-i) and $\gamma$ be a subprobability on $\FF$. There exists an estimator $\widehat s$ such that 
\begin{eqnarray*}
\lefteqn{C\E_{s}\cro{d^{2}(s,\widehat s)}}\quad\\
&\le& d^{2}(s,g\circ u)+\inf_{F\in\FF}\cro{d_{\infty}^{2}(g,F)+\tau(\Delta_{\gamma}(F)+
\D(F))}\\ &&\mbox{}+\tau\cro{\sum_{j=1}^{l}
k_{j}\log(1+2M_{j}R_{j}^{1/\beta_{j}})}+\sum_{j=1}^{l}\inf_{i\ge
1}\cro{l\w_{g,j}^{2}\pa{e^{-i}}+i\tau\pa{1+k_{j}\beta_{j}^{-1}}},
\end{eqnarray*}
for all $g\in\F_{l,c}$ and $u_{j}\in\gT_{j}$, $j=1,\ldots,l$.

In particular, for all $(\galpha,\gL)$-H\"olderian functions $g$ with $\galpha\in(0,1]^{l}$ and 
$\gL\in(\R_{+}^{\star})^{l}$, 
\begin{eqnarray}
C\E_{s}\cro{d^{2}(s,\widehat s)}&\le& d^{2}(s,g\circ u)+\inf_{F\in\FF}\cro{d_{\infty}^{2}(g,F)+
\tau(\Delta_{\gamma}(F)+\D(F))}\nonumber\\&&\mbox{}+\tau\sum_{j=1}^{l}
\cro{k_{j}\log(1+2M_{j}R_{j}^{1/\beta_{j}})+(\L_{j}\vee 1)\pa{1+k_{j}\beta_{j}^{-1}}},\qquad\quad 
\labe{Eq-F}
\end{eqnarray}
where  
\begin{equation}
\L_{j}={1\over 2\alpha_{j}}\log\pa{lL_{j}^{2}\over (1+k_{j}\beta_{j}^{-1})\tau}\quad\mbox{for }
j=1,\ldots,l.
\labe{Eq-Lj}
\end{equation}
\end{theorem}
Although this theorem is stated for a given value of $l$, we may, arguing as before, let $l$ vary 
and design a new estimator which achieves the same risk bounds (apart for the constant $C$) 
whatever the value of $l$.

As usual, the quantity $\inf_{F\in\FF}\cro{d_{\infty}^{2}(g,F)+\tau(\Delta_{\gamma}(F)+\D(F))}$ 
corresponds to the estimation rate for the function $g$ alone by using the collection $\FF$. In
particular, if $g\in\H^{\galpha}([-1,1]^{l})$ with $\galpha\in(\R_{+}^{\star})^{l}$, this bound is 
of order $\tau^{2\overline \alpha/(2\overline \alpha+l)}$ as $\tau$ tends to 0 for a classical choice
of $\FF$ (see Section~\ref{sect-ApproxS}). Since for all $j$, $g$ is also $(\alpha_{j}\wedge
1)$-H\"olderian as a function of $x_{j}$ alone, the last term in the right-hand side of~\eref{Eq-F},
which is of order $-\tau\log\tau$, becomes negligible as compared to $\tau^{2\overline
\alpha/(2\overline \alpha+l)}$ and therefore, when $s$ is really of the form $g\circ u$ with
$g\in\H^{\galpha}([-1,1]^{l})$ the rate we get for estimating $s$ is the same as that for
estimating $g$.\vspace{2mm}\\
%
{\em Proof of Theorem~\ref{T-para}:} 
For $\eta>0$ and $j=1,\ldots,l$, let $\Theta_{j}[\eta]$ be a maximal subset of $\Theta_{j}$ 
satisfying $\ab{t-t'}>\eta$ for all $t,t'$ in $\Theta_{j}[\eta]$. Since $\Theta_{j}$ is a subset of the
Euclidean ball in $\R^{k_{j}}$ centered at 0 with radius $M_{j}$, it follows from classical entropy
computations (see Lemma~4 in Birg\'e~\citeyearpar{MR2219712}) that $\log
\ab{\Theta_{j}[\eta]}\le k_{j}\log(1+2M_{j}\eta^{-1})$. For all $i\in\N^{\star}$, let $T_{j,i}$ be the
image of $\Theta_{j,i}=\Theta_{j}[(R_{j}e^{i})^{-1/\beta_{j}}]$ by the mapping $\gtheta\mapsto
u_{j}(\gtheta,.)$. Clearly, 
\[
\log \ab{T_{j,i}}\le \log\ab{\Theta_{j,i}}\le k_{j}\log\pa{1+2M_{j}R_{j}^{1/\beta_{j}}
e^{i/\beta_{j}}}\le k_{j}\cro{\log(1+2M_{j}R_{j}^{1/\beta_{j}})+i\beta_{j}^{-1}}
\]
and because of the maximality of $\Theta_{j,i}$ and~\eref{Eq-embed}, for all $\gtheta\in 
\Theta_{j}$ there exists $\overline \gtheta\in \Theta_{j,i}$ such that
$d(u_{j}(\gtheta,.),u_{j}(\overline \gtheta,.))\le R_{j}\ab{\gtheta-\overline \gtheta}^{\beta_{j}}\le
e^{-i}$ so that $T_{j,i}$ is an $e^{-i}$-net for $\gT_{j}$. For $j=1,\dots,l$, we set
$\TT_{j}=\bigcup_{i\ge 1}T_{j,i}$ so that the models in $\TT_{j}$  are merely the elements of the sets
$T_{j,i}$. For a model $T$ that belongs to 
$T_{j,i}\setminus\bigcup_{1\le i'<i}T_{j}[e^{-i'}]$ (with the convention $\bigcup_{\varnothing}=
\varnothing$) we set
\[
\Delta_{\lambda_{j}}(T)=\log\ab{T_{j,i}}+i\le  k_{j}\log\left(1+2M_{j}R_{j}^{1/\beta_{j}}\right)
+i\left(1+k_{j}\beta_{j}^{-1}\right)
\] 
which defines a measure $\lambda_{j}$ on $\TT_{j}$ satisfying
\[
\sum_{T\in \TT_{j}}\lambda_{j}(T)\le \sum_{i\ge 1}\sum_{t\in T_{j,i}}\lambda_{j}(\{t\})\le 
\sum_{i\ge 1}e^{-i}<1.
\]
Since for all $j$ and $T\in \TT_{j}$, $\D(T)=0$, we get the first risk bound by applying 
Theorem~\ref{T-general} to the corresponding set $\frak{S}$. To prove (\ref{Eq-F}), let us set 
$i(j)=\left\lfloor \L_{j} \right\rfloor\vee 1$ for $j=1,\ldots,l$ with $\L_j$ given by (\ref{Eq-Lj}), so 
that $1\le i(j)\le\L_{j}\vee 1$ and notice that, if $z\ge \L_{j}\vee 1$, then
$lL_{j}^{2}e^{-2\alpha_{j}z}\le z\tau\pa{1+k_{j}\beta_{j}^{-1}}$. If $\L_{j}\ge 1$, then $\L_{j}\le
i(j)+1\le 2\L_{j}$, hence
\[ 
lL_{j}^{2}e^{-2\alpha_{j}(i(j)+1)}\le (i(j)+1)\tau\pa{1+k_{j}\beta_{j}^{-1}}\le
2\L_{j}\tau \pa{1+k_{j}\beta_{j}^{-1}}
\]
and
\[
l\w_{g,j}^{2}(e^{-i(j)})\le lL_{j}^{2}e^{-2\alpha_{j}i(j)}\le 2e^{2\alpha_{j}}\L_{j}\tau \pa{1+k_{j}
\beta_{j}^{-1}}\le 2e^{2}\L_{j}\tau \pa{1+k_{j}\beta_{j}^{-1}}.
\]
Otherwise, $\L_{j}< 1$, $i(j)=1\ge \L_{j}\vee 1$ and 
\[
l\w_{g,j}^{2}(e^{-i(j)})\le lL_{j}^{2}e^{-2\alpha_{j}}\le \tau\pa{1+k_{j}\beta_{j}^{-1}},
\]
so that in both cases $l\w_{g,j}^{2}(e^{-i(j)})\le2e^{2}(\L_{j}\vee 1)\tau 
\pa{1+k_{j}\beta_{j}^{-1}}$, which leads to the conclusion.\cqfd

\subsubsection{Estimating a density by a mixture of Gaussian densities}
In this section, we consider the problem of estimating a bounded density $s$ with respect to 
some probability $\mu$ (to be specified later) on $E=\R^{k}$, $d$ denoting, as before, the 
$\IL_{2}$-distance on $\IL_{2}(E,\mu)$. We recall from Section~\ref{F2} that 
Theorem~\ref{T-models} applies to this situation with $\tau =n^{-1}\|s\|_{\infty}(1\vee\log
\|s\|_{\infty})$. A common way of modeling  a density on $E=\R^{k}$ is to assume that it is a
mixture of  Gaussian densities (or close enough to it). More precisely, we wish to approximate $s$ by
functions $\overline s$ of the form 
\begin{equation}
\overline s(x)=\sum_{j=1}^{l}q_{j}p(m_{j},\Sigma_{j},x)\ \ {\rm for\ all}\ \  x\in\R^{k},
\labe{Eq-MixG}
\end{equation}
where $l\ge 1$, $\gq=(q_{1},\ldots,q_{l})\in[0,1]^{l}$ satisfies $\sum_{j=1}^{l}q_{j}=1$ and for 
$j=1,\ldots,l$, $p(m_{j},\Sigma_{j},.)=d{\cal N}(m_{j},\Sigma_{j}^{2})/d\mu$ denotes
the density (with respect to $\mu$) of the Gaussian distribution ${\cal  N}(m_{j},\Sigma_{j}^2)$
centered at $m_{j}\in\R^{k}$ with covariance matrix $\Sigma_{j}^{2}$ for some symmetric 
positive definite matrix $\Sigma_{j}$. Throughout this section, we shall restrict to means $m_{j}$
with Euclidean norms not larger than some positive number $r$ and to matrices $\Sigma_{j}$ with
eigenvalues $\rho$ satisfying $\underline \rho\le \rho\le \overline \rho$ for positive numbers
$\underline \rho<\overline \rho$. In order to parametrize the corresponding densities, we
introduce the set $\Theta$ gathering the elements $\theta$ of the form $\theta=(m,\Sigma)$ where
$\Sigma$ is a positive symmetric matrix with eigenvalues in $[\underline \rho,\overline \rho]$
and $m\in\B_{k}(0,r)$. We shall consider $\Theta$ as a subset of $\R^{k(k+1)}$ endowed with the
Euclidean distance. In particular, the set $M_k$ of square $k\times k$ matrices of dimension $k$ is
identified to $\Bbb{R}^{k^2}$ and endowed with the Euclidean distance and the corresponding
norm $N$ defined by
\[
N^{2}(A)=\sum_{i=1}^{k}\sum_{j=1}^{k}A_{i,j}^2\quad\mbox{if }A=\left(A_{i,j}\st\right)
_{\stackrel{\scriptstyle{1\le i\le k}}{1\le j\le k}}.
\]
This norm derives from the inner product $[A,B]=\tr(AB^{*})$ (where $B^{*}$ denotes the
transpose of $B$) on $M_k$ and satisfies $N(AB)\le N(A)N(B)$ (by Cauchy-Schwarz inequality) 
and $N(A)=N(UAU^{-1})$ for all orthogonal matrices $U$. In particular, if $A$ is symmetric and
positive with eigenvalues bounded from above by $c$,  $N(A)\le \sqrt{k}c$. We shall use these 
properties later on. For $b=r^{2}/(2\overline \rho^{2})+k\log(\sqrt{2}\overline \rho/\underline
\rho)$ and $\mu$ the Gaussian distribution ${\mathcal N}(0,2\overline \rho^{2} I_{k})$ on
$\R^{k}$ (where $I_{k}$ denotes the identity matrix) we define the parametric set $\gT$ by
\[
\gT=\ac{u(\theta,.)=e^{-b/2}\sqrt{p(\theta,.)},\ \theta\in\Theta}.
\]
For parameters $\theta_{1}=(m_{1},\Sigma_{1}),\ldots,\theta_{l}=(m_{l},\Sigma_{l})$ in $\Theta$, 
the density $\overline s$ can be viewed as a composite function $g\circ u$ with 
\begin{equation}
g(y_{1},\ldots,y_{l})=e^{b}q_{1}y_{1}^{2}+\ldots+e^{b}q_{l}y_{l}^{2}
\labe{Eq-eqg}
\end{equation}
and $u=(u_{1},\ldots,u_{l})$ with $u_{j}(.)=u(\theta_{j},.)$ for $j=1,\ldots,l$. With our choices of $b
$ and $\mu$, $u(\theta,.)\in\T$ for all $\theta=(m,\Sigma)\in \Theta$ as required, since for all 
$x\in E$
\begin{eqnarray*}
p(\theta,x)&=&{(2\overline \rho^{2})^{k/2}\over \det \Sigma}\exp\cro{{\ab{x}^{2}\over 
4\overline \rho^{2}}-{\ab{\Sigma^{-1}(x-m)}^{2}\over 2}}\\ &\le&  (2\overline
\rho^{2}\underline \rho^{-2})^{k/2}\exp\cro{{\ab{x-m}^{2}\over 2\overline
\rho^{2}}+{\ab{m}^{2}\over 2\overline\rho^{2}}-{\ab{\Sigma^{-1}(x-m)}^{2}\over 2}}\\
&\le&  (2\overline \rho^{2}\underline \rho^{-2})^{k/2}e^{r^{2}/(2\overline \rho^{2})}\;\;\le\;\;e^{b}.
\end{eqnarray*}
An application of Theorem~\ref{T-para} leads to the following result. 
%
\begin{corollary}\lab{C-GM}
Let $s$ be a bounded density in $\IL_{2}(E,\mu)$, $\tau=n^{-1}\|s\|_{\infty}(1\vee\log
\|s\|_{\infty})$, $M=\sqrt{k}\overline \rho+r$, $b=r^{2}/(2\overline
\rho^{2})+k\log(\sqrt{2}\overline \rho/\underline \rho)$,
$R=\sqrt{k/2}e^{-b/2}\underline\rho^{-1}$ and 
\[
\L(\tau)={1\over 2}\log\pa{{4le^{2b}\tau^{-1}\over 1+k(k+1)}}.
\]
There exists an estimator $\widehat s$ satisfying for some universal constant $C>0$
\begin{equation}
C\E_{s}\cro{d^{2}(s,\widehat s)}\le \inf\cro{d^{2}(s,g\circ u)} + lk(k+1)\tau\cro{\log(1+2MR)
+(\L(\tau)\vee 1)},
\labe{Eq-para-final}
\end{equation}
where the infimum runs among all functions $u=(u_{1},\ldots,u_{l})\in~\gT^{l}$ and $g$ of the 
form~\eref{Eq-eqg}.
\end{corollary}
The second term in the right-hand side of~\eref{Eq-para-final} does not depend on $g$ nor $u$ and
is of order $-\tau\log \tau$ as $\tau$ tends to 0. As already mentioned, one can also consider 
many values of $l$ simultaneously and find the best one by using Theorem~\ref{T-melange}. Up to a
possibly different constant $C$, the risk of the resulting estimator then satisfies~\eref{Eq-para-final}
for all $l\ge 1$ simultaneously.  The problem of estimating the parameters involved in a mixture of
Gaussian densities in $\R^{k}$ has also been considered by Maugis and
Michel~\citeyearpar{Melange-Gauss}. Their approach is based on model selection among a family
of parametric models consisting of densities of the form~\eref{Eq-MixG}. Nevertheless, they restrict
to Gaussian densities with specific forms of covariance matrices only.\vspace{2mm}\\ 
\noindent{\em Proof of Corollary~\ref{C-GM}:} 
First note that for all $\theta\in \Theta$, $\ab{\theta}= \ab{m}+N(\Sigma)\le r+\sqrt{k}\overline \rho$. Hence, if  we can prove that for all $\theta_{0}=(m_{0},\Sigma_{0}),\theta_{1}=(m_{1},\Sigma_{1})$ in $\Theta$
\begin{equation}
d\pa{u(\theta_{0},.),u(\theta_{1},.)}\le {\sqrt{k/2}\ e^{-b/2}\over \underline \rho}\ab{\theta_{0}-\theta_{1}},
\labe{domi}
\end{equation}
Assumption~\ref{A-embed} will be satisfied with 
\[
M_{j}=M=r+\sqrt{k}\overline \rho\ \ \mbox{and}\ \ R_{j}=\sqrt{k/2}e^{-b/2}
\underline\rho^{-1}=R\quad\mbox{for } j=1,\ldots,l.
\]
We shall therefore be able to apply Theorem~\ref{T-para} with $\gT_{j}=\gT$ for all $j$, $\tau=
n^{-1}\|s\|_{\infty}(1\vee \log\|s\|_{\infty})$, $\FF=\ac{F}$ where $F$ is the linear span of 
dimension $\D(F)=l$ of functions $g$ of the form~\eref{Eq-eqg} and $\gamma$ the Dirac mass at
$F$. Since the functions $g$ of the form~\eref{Eq-eqg} are $\gL$-Lipschitz with $L_{j}=
2q_{j}e^{b}\le 2e^{b}$ for all $j$, we shall finally deduce~\eref{Eq-para-final} from~\eref{Eq-F}. We
therefore only have to prove~\eref{domi}. Let us first note that
\begin{equation}
d^{2}(u(\theta_{0},.),u(\theta_{1},.))=2e^{-b}h^{2}\pa{{\cal
N}(m_0,\Sigma_0^{2}),{\cal N}(m_1,\Sigma_1^{2})},
\labe{Eq-Lip}
\end{equation}
where $h$ denotes the Hellinger distance defined by (\ref{Eq-Hell}).
Some classical calculations show that 
\[
h^{2}\pa{{\mathcal N}(m_0,\Sigma_0^{2}),{\mathcal N}(m_1,\Sigma_1^{2})}=
1-{\exp\cro{-{1\over 4}\<m_{1}-m_{0},(\Sigma_{0}^{2}+\Sigma_{1}^{2})^{-1}(m_{1}-m_{0})\>}
\over \sqrt{\det\pa{{\Sigma_{0}^{-1}\Sigma_{1}+\Sigma_{0}\Sigma_{1}^{-1}\over 2}}}},
\]
and from the inequalities, $1-e^{-z}\le z$ and $\log(\det A)\le \tr (A-I_{k})$ which hold for all 
$z\in\R$ and all matrices $A$ such that $\det A>0$, by setting $\Sigma^{2}=\Sigma_{0}^{2}+
\Sigma_{1}^{2}$ we deduce that 
\begin{eqnarray*}
\lefteqn{4h^{2}\pa{{\mathcal N}(m_0,\Sigma_0^{2}),{\mathcal N}(m_1,\Sigma_1^{2})}}\quad\\
&\le& 2\log\cro{\det\pa{{\Sigma_{0}^{-1}\Sigma_{1}+\Sigma_{0}\Sigma_{1}^{-1}\over 2}}}+
\<m_{1}-m_{0},\Sigma^{-2}(m_{1}-m_{0})\>\\&\le&\tr\pa{\Sigma_{0}^{-1}\Sigma_{1}
+\Sigma_{0}\Sigma_{1}^{-1}-2I_{k}}+\<m_{1}-m_{0}, \Sigma^{-2}(m_{1}-m_{0})\>\\
&=&\tr\pa{(\Sigma_{0}-\Sigma_{1})\Sigma_{0}^{-1}(\Sigma_{0}-\Sigma_{1})\Sigma_{1}^{-1}}+
\<m_{1}-m_{0},\Sigma^{-2}(m_{1}-m_{0})\>\;\;=\;\;U_{1}+U_{2},
\end{eqnarray*}
with
\[
U_{1}=\tr\pa{(\Sigma_{0}-\Sigma_{1})\Sigma_{0}^{-1}(\Sigma_{0}-\Sigma_{1})\Sigma_{1}^{-1}}  \
\ {\rm and}\ \ U_{2}=\<m_{1}-m_{0},\Sigma^{-2}(m_{1}-m_{0})\>.
\]
It remains to bound $U_{1}$ and $U_{2}$ from above. For $U_{1}$, taking
$A=(\Sigma_{0}-\Sigma_{1})\Sigma_{0}^{-1}$ and $B=\Sigma_{1}^{-1}(\Sigma_{0}-\Sigma_{1})$
and using the fact that the eigenvalues of $\Sigma_{0}^{-1}$ and $\Sigma_{1}^{-1}$ are not larger
than $\underline \rho^{-1}$, we get
\begin{eqnarray*}
U_{1}&=&\cro{A,B}\;\;\le\;\;N(A)N(B)\;\;=\;\;N((\Sigma_{0}-\Sigma_{1})\Sigma_{0}^{-1})
N(\Sigma_{1}^{-1}(\Sigma_{0}-\Sigma_{1}))\\ &\le&
N(\Sigma_{0}^{-1})N(\Sigma_{1}^{-1})N^{2}(\Sigma_{0}-\Sigma_{1})\;\;\le\;\;
{kN^{2}(\Sigma_{0}-\Sigma_{1})\over \underline\rho^{2}}.
\end{eqnarray*}
Let us now turn to $U_{2}$. It follows from the same arguments that the symmetric matrix 
$\Sigma^{2}=\Sigma_{0}^{2}+\Sigma_{1}^{2}$ satisfies for all $x\in\R^{k}$, 
\[
\<\Sigma^{2} x,x\>=\ab{\Sigma_{0}x}^{2}+\ab{\Sigma_{1}x}^{2}\ge 2\underline\rho^{2}
\ab{x}^{2}, 
\]
hence
\[
U_{2}=\<m_{1}-m_{0},\Sigma^{-2}(m_{1}-m_{0})\>\le {\ab{m_{0}-m_{1}}^{2}\over 2\underline 
\rho^{2}}.
\]
Putting these bounds together, we obtain that 
\[
4h^2\pa{{\cal N}(m_0,\Sigma_0^{2}),{\cal N}(m_1,\Sigma_1^{2})}\le{k\over \underline 
\rho^{2}}\pa{N^{2}(\Sigma_{1}-\Sigma_{0})+\ab{m_{0}-m_{1}}^{2}}
={k\over \underline \rho^{2}}\ab{\gtheta_{0}-\gtheta_{1}}^{2},
\]
which, together with~\eref{Eq-Lip}, leads to~\eref{domi}.\cqfd

\Section{Proof of the main results\labs{H}}
Let us recall that, in this section, $d$ denotes the distance associated to the $\|\ \|_{q}$ norm of 
$\IL_{q}(E,\mu)$ and $d_{\infty}$ the distance associated to the supnorm on $\F_{l,\infty}$.

\subsection{Basic theorem\labs{H1}}
We shall first prove a general theorem of independent interest that applies to finite models
${\bf T}$ for functions in $\T^l$ and is at the core of all further developments. 
%
\begin{theorem}\lab{T-basic}
Let $I$ be a countable set and $\nu$ a subprobability on $I$. Assume that, for each
$\ell\in I$, we are given two countable families $\TT_\ell$ and $\FF_\ell$ of subsets of $\T^{l}$ and
$\F_{l,\infty}$ respectively such that each element ${\bf T}$ of $\TT_\ell$ is finite and each
$F\in\FF_\ell$ is a linear subspace of dimension $\D(F)\ge1$ of $\F_{l,\infty}$.  Let
$\lambda_\ell$ and $\gamma_\ell$ be subprobabilities on $\TT_\ell$ and $\FF_\ell$ respectively. One
can design an estimator $\widehat{s}=\widehat s(\bm{X})$  satisfying, for all $\ell\in I$, all
$u\in\T^{l}$ and $g\in\F_{l,c}$ with modulus of continuity $\w_g$, 
\begin{eqnarray*}
C\Bbb{E}_s\!\left[d^2\left(s,\widehat{s}\right)\st\right]&\le&\inf_{{\bf T}\in\TT_\ell}
\ac{l\inf_{t\in {\bf T}}\sum_{j=1}^l\w_{g,j}^2\left(\|u_j-t_j\|_p\right)+\tau
\cro{\Delta_{\lambda_\ell}({\bf T})+\log|{\bf T}|+\Delta_\nu(\ell)}}\nonumber\\&&\mbox{}
+d^2(s,g\circ u)
+\inf_{F\in\FF_\ell}\ac{d_{\infty}^{2}(g,F)+\tau\cro{\D(F)+\Delta_{\gamma_\ell}(F)}}.\qquad
\end{eqnarray*}
\end{theorem}
%
\noindent{\em Proof:}
For each $t\in\bigcup_{{\bf T}\in\TT_\ell} {\bf T}$ and $F\in\FF_\ell$ we consider the set
$F_t=\ac{f\circ t,\ f\in F}\subset\Bbb{L}_q(E,\mu)$, which is a ${\cal D}(F)$-dimensional
linear space. This leads to a new countable family of models $\SS_\ell$ together with a
subprobability $\pi_\ell$ on $\SS_\ell$ given by
\begin{equation}
\SS_\ell=\ac{F_{t},\,t\in \bigcup_{{\bf T}\in\TT_\ell} {\bf T},\,F\in\FF_\ell};\quad\pi_\ell(F_t)=
\gamma_\ell(F)\inf_{{\bf T}\in\TT_\ell\, ,\,{\bf T}\ni t} |{\bf T}|^{-1}\lambda_\ell({\bf T}).
\labe{Eq-SS}
\end{equation}
We then set
\[
\SS=\bigcup_{\ell\in I}\SS_\ell\qquad\mbox{and}\qquad\pi(F_t)=\nu(\ell)\pi_\ell(F_t)
\quad\mbox{for }F_t\in\SS_\ell.
\]
It follows that 
\[
\Delta_{\pi}(F_t)=\Delta_{\gamma_\ell}(F)+\inf_{{\bf T}\in\TT_\ell\, ,\,{\bf T}\ni t}
\left[\Delta_{\lambda_\ell}({\bf T})+\log(|{\bf T}|)\right]+\Delta_\nu(\ell)
\quad\mbox{for }F_t\in\SS_\ell.
\] 
Applying
Theorem~\ref{T-models} to $\SS$ and $\pi$ leads to an estimator $\widehat{s}$ satisfying, for
each $\ell\in I$,
\begin{eqnarray*}
\lefteqn{C\Bbb{E}_s\!\left[d^2\left(s,\widehat{s}\right)\st\right]}\\
&\le&\inf_{F\in\FF_\ell,\,{\bf T}\in\TT_\ell,\, t\in {\bf T}}\!\ac{d^2(s,F_t)+\tau\cro{\D(F)+
\Delta_{\gamma_\ell}(F)+\Delta_{\lambda_\ell}({\bf T})+\log|{\bf T}| +\Delta_\nu(\ell)}}.
\end{eqnarray*}
We now use Proposition~\ref{P-Approx} which implies that, for each $f\circ t\in F_{t}$,
\begin{eqnarray*}
d^{2}(s,f\circ t)&\le&\pa{\|s-g\circ u\|_q+\|g\circ u-f\circ t\|_q}^{2}\\ 
&\le&\pa{\|s-g\circ u\|_q+d_\infty(g,f)+2^{1/q}\sum_{j=1}^{l}
\w_{g,j}\left(\|u_{j}-t_{j}\|_q\right)}^{2}\\&\le&3\pa{\|s-g\circ u\|_q^{2}
+d^2_\infty(g,f)+4l\sum_{j=1}^{l}\w_{g,j}^{2}\left(\|u_{j}-t_{j}\|_q\right)},
\end{eqnarray*}
for some universal constant $C$ since $2^{1/q}\le2$. The conclusion follows from a
minimization over all  possible choices for $f$ and $t$.\cqfd
%

\subsection{Building new models\labs{H2}}
In order to use Theorem~\ref{T-basic}, which applies to finite sets ${\bf T}$, starting from the
models $T$ which satisfy Assumption~\ref{A-models}, we need to derive new models from the
original ones. Let us first observe that, since $u_j$ takes its values in $[-1,1]$ and $\mu$ is a
probability on $E$, $d(0,u_j)\le1$. It is consequently useless to try to approximate $u_j$ by
elements of $\Bbb{L}_q(E,\mu)$ that do not belong to ${\cal B}(0,2)$ since 0 always does better. 
We may therefore replace $T\subset\mathbb{L}_q(E,\mu)$ by $\left(T\cap{\cal
B}(0,2)\st\right)\cup\{0\}$, denoting again the resulting set, which remains a subset of some
$\D(T)$-dimensional linear space, by $T$. Moreover, this modification can only decrease the
value of $d(T,u_j)$. Since now
$T\subset{\cal B}(0,2)$, we can use the discretization argument described by the following
lemma.
%
\begin{lemma}\lab{L-approx1}
Let $T\subset{\cal B}(0,2)$ be either a singleton (in which case $\D(T)=0$) or a subset of some
$\D(T)$-dimensional linear subspace of $\Bbb{L}_q(E,\mu)$ with $\D(T)\ge1$. For each
$\eta\in(0,1]$, one can find a subset $T[\eta]$ of $\T$ with cardinality bounded by
$(5/\eta)^{\D(T)}$ such that
\begin{equation}
\inf_{t\in T[\eta]}d(t,v)\le\inf_{t\in T}d(t,v)+[\eta\wedge\D(T)]\quad\mbox{for all }v\in\T.
\labe{Eq-approx1}
\end{equation}
\end{lemma}
%
\noindent{\em Proof:}
If $\D(T)=0$, then $T=\{t\}$, we set $T[\eta]=\{(-1\vee t)\wedge1\}$ and the result is
immediate since $v$ takes its values in $[-1,1]$. Otherwise, let $T'$ be a maximal subset of $T$
such that $d(t,t')>\eta$ for each pair $(t,t')$ of distinct points in $T'$. Then, for each $t\in T$
there exists $t'\in T'$ such that $d(t,t')\le\eta$ and it follows from Lemma~4 in Birg\'e~\citeyearpar{MR2219712}
that $|T'|\le(5/\eta)^{\D(T)}$. Now set $T[\eta]=\{(-1\vee t)\wedge1,\,  t\in T'\}$.
Then (\ref{Eq-approx1}) holds since $\D(T)\ge1$. \cqfd\vspace{2mm}\\
We are now in a position to build discrete models for approximating the elements of $\T^l$.
Given $j\in\ac{1,\ldots,l}$, $T_j$ in $ \TT_{j}$ and some $i\in\N^\star$, the previous lemma 
provides a set $T_j\!\left[e^{-i}\right]$ satisfying $\left|T_j\!\left[e^{-i}\right]\right|\le
\exp\cro{\D(T_j)\pa{i+\log 5}}$. Moreover, 
\begin{equation}
d\left(u_{j},T_j\left[e^{-i}\right]\right)\le d\left(u_j,T_j\right)+\left[e^{-i}\wedge
{\cal D}(T_j)\right]\quad\mbox{for all }u\in\T^l\mbox{ and }i\in\N^\star.
\labe{Eq-Newapp}
\end{equation}
We then define the family $\TT$ of models by
\begin{equation}
\TT=\ac{{\bf T}=\prod_{j=1}^lT_j\!\left[e^{-i_j}\right]
\mbox{ with }(i_j,T_j)\in\Bbb{N}^\star\times\TT_j\mbox{ for }j=1,\ldots,l}.
\labe{Eq-def-T}
\end{equation}
Then each ${\bf T}=T_{1}\left[e^{-i_1}\right]\times \ldots \times T_l\left[e^{-i_l}\right]$ in 
$\TT$ has a finite cardinality bounded by
\begin{equation}
\log\ab{{\bf T}}\le\sum_{j=1}^l\D(T_j)\pa{i_j+\log 5}.
\labe{Eq-maj-logT}
\end{equation}

\subsection{Proof of Theorem~\ref{T-general}\labs{H3}}
Starting from the families $\TT_j$, $1\le j\le l$, we build the set $\TT$ given by
(\ref{Eq-def-T}) as indicated in the previous section and we apply Theorem~\ref{T-basic} to
$\FF$ and $\TT$. This requires to define a suitable subprobability $\lambda$ on $\TT$, which
can be done by setting, for each ${\bf T}=T_{1}\left[e^{-i_1}\right]\times
\ldots \times T_l\left[e^{-i_l}\right]$ in $\TT$,
\[
\lambda({\bf T})=\prod_{j=1}^{l}\lambda_j(T_j)\exp\left[-i_j{\cal D}(T_j)\right]
\quad\mbox{or}\quad
\Delta_\lambda({\bf T})=\sum_{j=1}^l\left[\Delta_{\lambda_j}(T_j)+i_j{\cal D}(T_j)\right].
\]
Applying Theorem~\ref{T-basic} to $\FF$ and $\TT$ with $ I$ reduced to a single element
and $\nu$ the Dirac measure and using (\ref{Eq-maj-logT}) and (\ref{Eq-Newapp}) which
implies that
\begin{eqnarray*}
\inf_{t_j\in T_j\left[e^{-i_j}\right]}\w_{g,j}\left(\|u_j-t_j\|_p\right)&\le&
\w_{g,j}\left(\left[e^{-i_j}\wedge{\cal D}(T_j)\right]+d\left(u_j,T_j\right)\right)\\&\le&
\w_{g,j}\left(e^{-i_j}\wedge{\cal D}(T_j)\right)+\w_{g,j}\left(d(u_j,T_j)\st\right)
\end{eqnarray*}
by the subadditivity property of the modulus of continuity $\w_{g,j}$, we get the risk bound
\begin{eqnarray*}
C\Bbb{E}_s\!\left[d^2\left(s,\widehat s)\right)\st\right]&\le&
\inf\sum_{j=1}^l\left\{\rule{0mm}{4mm}2l\left[\w^2_{g,j}\!\left(d(u_j,T_{j})\right)+
\w^2_{g,j}\!\left(e^{-i_j}\wedge{\cal D}(T_j)\st\right)\right]\right.\\&&\left.
\qquad\quad\,\mbox{ }+\tau\left[\Delta_{\lambda_{j}}(T_{j})+
(2i_j+\log5)\D(T_{j})\right]\rule{0mm}{4mm}\right\}\\&&\mbox{}+d^2(s,g\circ u)
+\inf_{F\in\FF}\ac{d_{\infty}^{2}(g,F)+\tau\cro{\D(F)+\Delta_{\gamma}(F)}},
\end{eqnarray*}
where the first infimum runs among all  $T_j\in\TT_{j}$ and all $i_j\in\N^\star$ for
$j=1,\ldots,l$. Setting $i_j=i(g,j,T_j)$ implies that $l\w^2_{g,j}\left(e^{-i_j}\wedge
{\cal D}(T_j)\right)\le\tau i_j{\cal D}(T_j)$, which proves (\ref{Eq-MainRB}). As to
(\ref{Eq-RB2a}), it simply derives from the fact that, if ${\cal D}(T)\ge1$, then
\[
i(g,j,T)\le \left\lceil(2\alpha_j)^{-1}\log\left(lL_j^2[\tau\D(T)]^{-1}\right)\right\rceil
\le\left[\alpha_j^{-1}\log\left(lL_j^2[\tau\D(T)]^{-1}\right)\right]\bigvee1=\L_{j,T}.
\]
%

\subsection{Proof of Theorem~\ref{T-melange}\labs{H4}}
It follows exactly the line of proof of Theorem~\ref{T-general} via Theorem~\ref{T-basic} with
an additional step in order to mix the different families of models corresponding to the various
sets $\frak{S}_{\ell}$. To each $\frak{S}_{\ell}$ corresponds a family of models $\SS_{\ell}$ and a
subprobability $\pi_{\ell}$ on $\SS_{\ell}$ given by (\ref{Eq-SS}). We again apply 
Theorem~\ref{T-basic} with $ I$ and $\nu$ as given in
Theorem~\ref{T-melange}.

\subsection{Proof of Lemma~\ref{L-opt}\labs{H5}}
If $D=1$, we get the bound $a+b$. When $a>b$, we can choose $D$ such that
$(a/b)^{1/(\theta+1)}\le D<(a/b)^{1/(\theta+1)}+1$, so that
\[
aD^{-\theta}+bD<a(a/b)^{-\theta/(\theta+1)}+b\left[(a/b)^{1/(\theta+1)}+1\right]
=b+2a^{1/(\theta+1)}b^{\theta/(\theta+1)}
\]
and the bound $b+\left[2a^{1/(\theta+1)}b^{\theta/(\theta+1)}\wedge a\right]$ follows. If
$b\ge a$, the bound $2b$ holds, otherwise $b<a^{1/(\theta+1)}b^{\theta/(\theta+1)}$ and the
conclusion follows.

\subsection{Proof of Proposition~\ref{P-Approx}\labs{H6}}
It relies on the following lemma the proof of which is postponed to the end of the section. 
%
\begin{lemma}\lab{L-approx}
Let $(E,{\cal E},\mu)$ be some probability space and $w$ some nondecreasing and
nonnegative concave function on $\Bbb{R}_+$ such that $w(0)=0$. For all
$p\in[1,+\infty]$ and $h\in\Bbb{L}_p(\mu)$, 
\[
\|w(|h|)\|_p\le 2^{1/p}w(\|h\|_p),
\]
with the convention $2^{1/\infty}=1$.
\end{lemma}
We argue as follows. For all $y,y'\in [-1,1]^l$, $|g(y)-g(y')|\le\sum_{j=1}^l
\w_{g,j}\left(\left|y_{j}-y_{j}'\right|\right)$ and, since $\mu$ is a probability on $E$, 
\begin{eqnarray*}
\|g\circ u-f\circ t\|_p&\le&\|g\circ u-g\circ t\|_p+\|g\circ t-f\circ t\|_p\\&\le&
\left\|\sum_{j=1}^l\w_{g,j}(|u_{j}-t_{j}|)\right\|_p+\|g\circ t-f\circ t\|_p\\&\le&
\sum_{j=1}^l\|\w_{g,j}(|u_{j}-t_{j}|)\|_p+\sup_{y\in[-1,1]^l}|g(y)-f(y)|\\&\le& 
2^{1/p}\sum_{j=1}^l\w_{g,j}\left(\|u_{j}-t_{j}\|_p\right)+d_{\infty}(g,f),
\end{eqnarray*}
which proves the proposition. Let us now turn to the proof of the lemma.\vspace{2mm}\\
\noindent{\em Proof of Lemma~\ref{L-approx}:}
Since there is nothing to prove if $\|h\|_p=0$, we shall assume that $\|h\|_p>0$. The
assumptions on $w$ imply that, for all $0<a<b$, $b^{-1}w(b)\le a^{-1}w(a)$ and $w(a)\le
w(b)$. Consequently, for $p\in [1,+\infty[$, 
\begin{eqnarray*}
\int_E w^p(|h|)\,d\mu &=& \int_E w^p(|h|)\1_{|h|\le b}\,d\mu+
\int_E w^p(|h|)\1_{|h|>b}\,d\mu\\ &\le& w^p(b)+
\int_E\frac{w^p(|h|)}{|h|^p}|h|^p\1_{|h|> b}\,d\mu\;\;\le\;\;
w^p(b) +\frac{w^p(b)}{b^p}\int_E|h|^p\,d\mu,
\end{eqnarray*}
and the result follows by choosing $b=\|h\|_p$. The case $p=\infty$ can be deduced
by letting $p$ go to $+\infty$.\cqfd

\subsection{Proof of Proposition~\ref{P-Holder}\labs{H7}}
It suffices to show that for all $i\in\ac{1,\ldots,k}$ and $x\in[-1,1]^{k}$, the map  $g\circ
u_{x}(t)=g\circ u(x_{1},\ldots,x_{i-1},t,x_{i+1},\ldots,x_{k})$ from $[-1,1]$ into $\R$ belongs to
${\cal H}^{\gtheta}([-1,1]^{k})$. If at least $\alpha$ or $\beta_{i}$ are not larger than 1, the
result is clear. Otherwise both are larger than 1 and we can write $\beta_{i}=b_{i}+\beta_{i}'$ and
$\alpha=a+ \alpha'$ with $a,b\in\N^\star$ and $\beta_{i}', \alpha'\in (0,1]$. Both functions 
$g$ and $u_{x}$ are   $b_{i}\wedge a$ times differentiable and the derivatives
$g^{(\ell)}\circ u_{x}$ and $u_{x}^{(\ell)}$ for $\ell=0,\ldots, b_{i}\wedge a$ are H\"olderian with
smoothness $\rho=(\beta_{i}-b_{i}\wedge \alpha)\wedge (\alpha- b_{i}\wedge a)\in(0,1]$. Since
the derivative of order  $b_{i}\wedge a$ of $g\circ u_{x}$ is a polynomial with respect to these
functions, we derive (\ref{Eq-compsmooth}) from the fact that the set $(\H^{\rho}([-1,1]^{k}),+,.)$
is an algebra on $\R$. 

We shall prove the second part of the proposition for the case $k=1$ only since the general case can be proved by similar arguments. For $\rho>0$, let $h_{\rho}\in\H^{\rho}([-1,1])\setminus \bigcup_{\rho'>\rho}\H^{\rho'}([-1,1])$. Given $\alpha,\beta>0$, we distinguish between the cases below and the reader can check that for each of these  $g\in \H^{\alpha}([-1,1])$, $u\in \H^{\beta}([-1,1])$, $g\circ u\in \H^{\theta}([-1,1])$ with $\theta=\phi(\alpha,\beta)$ but $g\circ u\not \in \H^{\theta'}([-1,1])$ whatever $\theta'>\theta$. If $\alpha,\beta\le 1$, take $g(x)=|x|^{\alpha}$ and $u(y)=|y|^{\beta}$ for all $x,y\in[-1,1]$, if $1<\beta$ and $\alpha\le \beta$, take $g=h_{\alpha}$ and $u(y)=y$ for all $y\in[-1,1]$, finally, if $\alpha>1$ and $\alpha> \beta$, take $g(x)=x$ for all $x\in [-1,1]$ and $u=h_{\beta}$.

\subsection{Proof of Lemma~\ref{L-comphold}\labs{H8}}
For all $\alpha>0$, the map defined for $y$ in $(0,+\infty)$ by
\[
\phi_{\alpha}(y)={1\over\phi(\alpha,1/y)}=
\left\{\begin{array}{ll}y(\alpha\wedge 1)^{-1}&\mbox{if }y\ge (\alpha\vee 1)^{-1};\\
\alpha^{-1}&\mbox{otherwise,}\end{array}\right.
\]
is positive, piecewise linear and convex. Hence, 
\[
{1\over\overline\theta}={1\over k}\sum_{i=1}^{k}\phi_{\alpha}\left(\frac{1}{\beta_i}\right)
\ge\phi_\alpha\left(\frac{1}{\overline\beta}\right)
={1\over\phi(\alpha,\overline\beta)}
\]
and equality holds if and only if $\beta_{i}\le (\alpha\vee 1)$ for all $i$ or if for all $i$, 
$\beta_{i}\ge (\alpha\vee 1)$.  We conclude by using the fact that $\phi(\alpha,z)\le 
z(\alpha\wedge 1)$ for all positive number $z$ and that equality holds if and only if $z\le
\alpha\vee 1$.

\bibliographystyle{apalike}

\end{document}